\documentclass[preprint,12pt]{elsarticle}
\usepackage[hidelinks]{hyperref}
\usepackage{geometry}
\usepackage{multirow}
\usepackage{array}
\usepackage{float}
\usepackage{amsmath,amssymb,amsthm,mathtools,bm}
\usepackage{graphicx}
\usepackage{booktabs}
\usepackage{enumitem}
\usepackage{tikz}
\usetikzlibrary{shapes.geometric, arrows.meta, positioning, calc, fit}

\begin{document}

\begin{frontmatter}

\title{Efficient and Accurate Model Order Reduction for Integral Electromagnetic Formulations in Fusion Device Transient Analysis Toward AI-Enabled Modeling}

\author[inst1]{Salvatore Ventre}

\affiliation[inst1]{
    organization={Department of Electrical and Information Engineering, University of Cassino and Southern Lazio},
    addressline={Via Gaetano di Biasio 43},
    city={Cassino},
    postcode={03043},
    country={Italy}
}

\begin{abstract}
The numerical simulation of electromagnetic transients in fusion devices is essential for analyzing plasma stability and disruptive events. However, it remains computationally demanding due to the large-scale dense systems arising from integral formulations. 

This work proposes a model order reduction (MOR) strategy for transient electromagnetic problems based on integral formulations. Unlike operator-based compression techniques (such as $\mathcal{H}$-matrix approaches), the reduced space is constructed directly from the transient excitation. In contrast to classical snapshot- and transfer-function-based MOR approaches, the proposed formulation avoids repeated explicit inversions or factorizations of the dense integral operator during the MOR basis-construction stage. By combining wavelet-based temporal compression with source-driven Krylov projections, the method generates reduced models tailored to the dynamically reachable responses of the prescribed excitation families.

Numerical validations on various plasma events and fusion-relevant scenarios demonstrate the robustness of the strategy, achieving substantial computational speedups while accurately preserving the transient electromagnetic response. Finally, the method is successfully applied to the null-field problem to efficiently generate training data for neural-network surrogates, contributing toward physics-consistent AI-enabled fusion modeling.
\end{abstract}



\end{frontmatter}

\section{Introduction}
\label{sec:introduction}

In fusion technology, the numerical simulation of electromagnetic phenomena is essential for analyzing plasma stability and predicting transient events in devices such as ITER~\cite{iter_project}. Fusion systems are inherently complex due to the strong coupling between plasma dynamics and surrounding conducting or magnetic structures. From an engineering perspective, this complexity strongly impacts both structural design and the development of robust control strategies.

Among the most critical plasma events are disruptions—such as vertical displacement events (VDEs) and asymmetric VDEs (AVDEs)~\cite{hender2007}. These events induce intense eddy currents and severe electromagnetic loads, threatening the structural integrity of the device. Accurately simulating these fast three-dimensional transient phenomena requires high-fidelity numerical models. However, the resulting massive spatial discretizations and small time steps make physics-based simulations computationally prohibitive for real-time control, optimization, and digital-twin frameworks~\cite{ferraro2020advancing}.

To address these challenges, integral formulations of the magneto--quasi--static (MQS) problem, such as the CARIDDI code~\cite{albanese1988,bossavit1981,passarotto2022}, provide an efficient computational framework. By restricting the discretization to conducting regions through an electric vector potential formulation, these methods avoid meshing large vacuum domains. Despite these advantages, the resulting discretized integral systems remain large-scale and dense, scaling cubically during matrix factorization and quadratically during transient time stepping.

Several matrix-compression and acceleration techniques—including Fast Multipole Methods (FMM), hierarchical (\(\mathcal{H}\)- and \(\mathcal{H}^{2}\)-) matrices, low-rank approximations, and symmetry-based reductions such as the \(n\)-fold rotational formulation—have been developed to handle large-scale dense integral operators~\cite{greengard1987fast,borm2003,rubinacci2004fast,ma2018,rubinacci2020application,ventre2022hmatrix}. While these methods efficiently reduce memory occupation and accelerate matrix--vector operations associated with the dense inductance operator~\cite{rubinacci2002,rubinacci2008fast}, they remain essentially operator-driven. Since the compression acts on the algebraic structure of the integral operator independently of the transient excitation manifold, the reduction of the dynamically relevant solution space remains intrinsically limited. As a result, for realistic fusion-device applications, the achievable computational speedups often remain moderate, with more substantial gains emerging only for extremely large-scale problems.

Recently, artificial intelligence (AI) has emerged as a promising tool to overcome these computational limitations in the design, simulation, and real-time monitoring of multiphysics fusion systems. However, its application in safety-critical reactors is limited by the black-box nature of purely data-driven approaches, which often struggle with interpretability, limited training data, and highly nonlinear dynamics~\cite{ferraro2020advancing}. Consequently, future AI frameworks must remain tightly integrated with physics-based predictive architectures in order to preserve consistency with the underlying plasma and electromagnetic dynamics.

To achieve this integration, model order reduction (MOR)~\cite{antoulas2005approximation,benner2015survey} provides a powerful framework by reducing the dimension of the governing system while preserving its dominant physical properties. By enabling fast and reliable simulations, MOR naturally acts as a bridge between high-fidelity physics-based models and the computational efficiency required by digital twins. In the present work, this concept is extended through an excitation-driven MOR strategy, where the reduced space is constructed to capture the transient dynamics effectively reachable under the considered source families. In this way, the reduction process is directly tailored to the physically relevant solution manifold rather than to the full operator, thereby ensuring highly effective compression. Due to the significantly reduced dimensionality of the resulting MOR strategy, the method can be efficiently exploited to generate training datasets for neural-network surrogates.

The remainder of the paper is organized as follows. Section~\ref{sec:mor_soa} reviews the state of the art in model order reduction for transient electromagnetic problems in fusion devices. Section~\ref{sec:full_order_nullspace} introduces the mathematical formulation of the full-order integral electromagnetic problem. Section~\ref{sec:null_field_application} presents the proposed excitation-driven Krylov MOR strategy and details its application to the null-field control problem, including the POD--NN surrogate architecture. Section~\ref{sec:numerical_results} reports the numerical experiments and discusses the corresponding validation results. Finally, Section~\ref{sec:conclusion} summarizes the main conclusions and outlines future research directions.

\section{Model Order Reduction: state of the art}
\label{sec:mor_soa}

Model order reduction (MOR) aims to construct low-dimensional models that accurately reproduce the input--output behavior of large-scale dynamical systems at a fraction of the computational cost~\cite{antoulas2005approximation,benner2021mor}. In computational electromagnetics, MOR techniques are broadly divided into two main categories: system-theoretic and snapshot-based projection methods~\cite{benner2021mor}.

Among system-theoretic approaches, classical moment-matching and Krylov-subspace methods~\cite{golub2013matrix} construct reduced spaces in the frequency domain by matching the moments of the transfer function around selected expansion points~\cite{freund2003model,grimme1997krylov}. In electromagnetic applications, these techniques are typically implemented through Arnoldi-type algorithms~\cite{bonotto2017}.

Conversely, snapshot-based methods—such as Proper Orthogonal Decomposition (POD) and Reduced Basis approaches—extract the projection space from representative full-order transient solutions. In fusion applications, POD techniques are increasingly adopted to simulate plasma instabilities~\cite{lucchini2024mor}, capturing the system dynamics directly in the time domain under selected operating conditions.

When applied to large-scale integral formulations, however, both classes share a critical computational bottleneck: they inherently require repeated inversions or factorizations of dense integral operators. Whether through frequency-domain expansions or time-domain snapshot generation, handling such large dense matrices rapidly becomes computationally prohibitive for ITER-scale configurations.

These limitations become particularly severe in fusion applications, where plasma-driven electromagnetic transients exhibit strong temporal variability and rapid scenario evaluation is essential. To address these challenges, this work proposes an excitation-driven MOR strategy, where the reduced space is constructed directly from a compressed representation of the transient forcing manifold.

Unlike classical transfer-function-based or snapshot-based approaches, the proposed formulation constructs the reduced space without requiring repeated dense factorizations of the integral operator during basis generation, while also eliminating the need for preliminary full-order transient simulations during the basis-construction stage. As shown throughout the paper, the methodology combines a compact representation of the RHS with a Krylov approximation based on sequential applications of \(R^{-1}L\), where \(R\) and \(L\) denote the sparse resistance and dense inductance operators of the full-order model, respectively (see Appendix~\ref{app:cariddi}). In this way, the method generates source-informed reduced spaces tailored to the dominant transient responses associated with the prescribed plasma events.

\section{Passive structures electromagnetic model and null-space formulation}
\label{sec:full_order_nullspace}

The ultimate objective is to simulate the fully self-consistent problem, where plasma dynamics and surrounding conducting structures interact through mutual electromagnetic coupling~\cite{holzl2012,isernia2023}. However, due to its intrinsic complexity, the present work adopts a non-self-consistent approach with a prescribed plasma evolution. This allows the analysis to focus entirely on the fundamental characteristics of the proposed MOR strategy and to assess its efficiency before addressing the fully coupled self-consistent problem. For completeness, a brief overview of the underlying eddy-current integral formulation is provided in Appendix~\ref{app:cariddi}.

The following differential-algebraic system governs the electromagnetic transient problem:
\begin{align}
R\,I(t) + L\,\frac{dI}{dt}(t) + F^T\Phi(t) &= b_i(t), \label{eq:full1}\\
F\,I(t) &= \alpha_{J_0}(t), \label{eq:full2}
\end{align}

where \(I(t)\) is the vector of induced currents, representing the numerical counterpart of the induced current density \(\mathbf{J}(t)\), \(\Phi(t)\) contains the unknown electrode voltages associated with current-driven electrodes, \(R\) is the sparse resistive operator, \(L\) is the dense inductive operator, and \(F\) is a  constraint operator. Both \(R\) and \(L\) are symmetric positive-definite matrices of size \(N\), where \(N\) denotes the total number of degrees of freedom associated with the spatial discretization mesh.


The forcing term \(\alpha_{J_0}(t)\) represents the imposed currents directly enforced through the algebraic constraint equation~\eqref{eq:full2}, while the remaining electromagnetic excitations acting on the passive conducting structures are collected in \(b_i(t)\), namely
\begin{equation}
b_i(t)
=
V_{\mathrm{volt}} \alpha_{\mathrm{volt}}(t)
+
V_{\mathrm{3D}} \alpha_{\mathrm{3D}}(t)
+
V_{\mathrm{axi}} \alpha_{\mathrm{axi}}(t),
\label{eq:forcing}
\end{equation}
where the three contributions correspond to voltage-driven electrodes, non-axisymmetric 3D conductors, and axisymmetric sources, respectively. In particular, the axisymmetric excitations include poloidal-field coils, central-solenoid coils, and equivalent plasma current filaments, while the non-axisymmetric sources arise from external 3D conductors such as TF coils. This formulation allows the representation of plasma--wall current exchange phenomena, including HALO current effects~\cite{bettini2013}. In summary, four families of electromagnetic excitations are considered: axisymmetric sources, non-axisymmetric 3D sources, voltage-driven electrodes, and imposed current sources.

This decomposition separates spatial and temporal contributions, allowing the geometrical operators to be assembled only once and reused across different transient scenarios by updating only the excitation waveforms.

To eliminate the constraint matrix \(F\) and the unknown voltages \(\Phi(t)\) from~\eqref{eq:full1}, the current vector is decomposed into null-space coordinates:
\begin{equation}
I(t) = K\,y(t) + I_0(t),
\qquad
\text{subject to } FK = 0,
\label{eq:null_decomp_continuous}
\end{equation}
where \(K\) spans the right null-space of \(F\), \(y(t)\) represents the unconstrained state coordinates, and \(I_0(t)\) is a particular solution satisfying~\eqref{eq:full2}. Substituting~\eqref{eq:null_decomp_continuous} into~\eqref{eq:full1} and projecting the system onto the range of \(K^T\) gives
\begin{equation}
R_K\,y(t)
+
L_K\,\frac{dy}{dt}(t)
=
K^T b_i(t)
-
K^T R\,I_0(t)
-
K^T L\,\frac{d I_0}{dt}(t),
\label{eq:nullspace_continuous}
\end{equation}
with
\[
R_K = K^T R K,
\qquad
L_K = K^T L K.
\]

This transformation removes the saddle-point structure of the original algebraic system, yielding symmetric positive-definite projected operators \(R_K\) and \(L_K\), and thereby improving both numerical stability and computational efficiency.

When electrodes are absent, the standard unconstrained formulation is naturally recovered by setting
\[
K = I,
\qquad
I_0(t)=0.
\]

For realistic large-scale applications such as ITER, solving~\eqref{eq:nullspace_continuous} remains computationally intensive. All stages of the CARIDDI code—including matrix assembly, factorization, transient solution, RHS evaluations, and postprocessing—are implemented in parallel using distributed-memory HPC architectures through MPI~\cite{mpi_standard} and ScaLAPACK~\cite{scalapack}. As discussed in Appendix~\ref{app:transient}, the transient solver can exploit either a direct approach based on matrix factorization~\cite{fresa2005eddy} or a modal approach based on generalized eigenvalue decomposition~\cite{ventre2025modal}.

The total computational cost of the transient simulation can be decomposed into two main contributions:
\[
T_{\mathrm{trans}}
=
T_{\mathrm{setup}}
+
T_{\mathrm{march}},
\]
where \(T_{\mathrm{setup}}\) denotes the preprocessing stage, performed only once before the transient analysis, while \(T_{\mathrm{march}}\) represents the cumulative cost associated with the time-marching procedure over all time steps.

The setup phase includes operations such as matrix factorizations or generalized eigenvalue decompositions and can be expressed as
\[
T_{\mathrm{setup}}
=
K_s N_{\mathrm{full}}^3,
\]
where the proportionality constant \(K_s\) depends on the adopted strategy and is typically much smaller for matrix factorizations than for generalized eigensolutions~\cite{fresa2005eddy,golub2013matrix,ventre2025modal}.

The marching cost depends on the number of time steps \(N_t\):
\[
T_{\mathrm{march}}
=
N_t\,T_{\mathrm{step}}.
\]

For a direct time-marching scheme and for a modal formulation, the per-step costs can be written as
\[
T_{\mathrm{step}}^{\mathrm{direct}}
=
K_d N_{\mathrm{full}}^2,
\qquad
T_{\mathrm{step}}^{\mathrm{modal}}
=
K_m N_{\mathrm{full}},
\]
respectively, where \(K_d\) and \(K_m\) are implementation-dependent constants.

When the number of time steps \(N_t\) becomes sufficiently large --- which is typically the case in practical fusion applications characterized by long transients, large final times, and small time steps --- the one-time setup cost becomes comparatively negligible with respect to the cumulative marching contribution. In this regime, the modal approach becomes particularly advantageous, since its marching cost grows linearly with the system size, whereas the direct formulation scales quadratically.

Thus, despite the higher setup cost associated with the eigensolution, the modal formulation can provide a significant computational advantage for long transients~\cite{ventre2025modal}.

\section{Introduction to the excitation-driven Krylov MOR}
\label{sec:intro_krylov_mor}

Traditional MOR techniques for electromagnetic integral formulations are predominantly based on the approximation of the system operator. Matrix-compression approaches, such as hierarchical matrices and related low-rank methods~\cite{greengard1987fast,ventre2022hmatrix,borm2003,cheng1999fast,ma2018}, approximate the dense integral operator independently of the excitation sources. Similarly, classical Krylov and moment-matching techniques~\cite{grimme1997krylov,bonotto2017,freund2003model} generate the reduced basis directly from transfer-function expansions, requiring repeated applications of the dynamical operator inverse. In contrast, the present formulation adopts an excitation-driven strategy, where the reduced space is constructed directly from the transient forcing manifold associated with the external sources.

The proposed MOR framework relies on two main pillars: (i) a Krylov-based approximation of the transient electromagnetic dynamics, and (ii) a reduced representation of the excitation space. The first component captures the dominant transient dynamics of the magneto--quasi--static (MQS) system, while the second exploits the strong redundancy typically present in transient source excitations of fusion devices.

The rationale behind the Krylov subspace construction can be understood from the transient evolution of a simple RL circuit governed by the following dynamic equation 
\[
L\frac{dI}{dt}+RI(t)=e(t),
\]
where \(L\) and \(R\) denote the inductance and resistance matrices, respectively. Applying a backward Euler discretization with time step \(\Delta t\) and introducing the current increment
\[
\Delta I_n = I_n-I_{n-1},
\]
yields
\[
(L+\Delta t\,R)\,\Delta I_n
=
\Delta t\,e_n
-
\Delta t\,R\,I_{n-1}.
\]

Multiplying by \(L^{-1}\) and defining the inverse time-constant operator
\[
\tau^{-1}=L^{-1}R,
\]
the algebraic system becomes
\[
\left(1+\Delta t\,\tau^{-1}\right)\Delta I_n
=
\Delta t\,L^{-1}e_n
-
\Delta t\,\tau^{-1}I_{n-1}.
\]

Introducing the effective forcing term
\[
b_n=L^{-1}e_n-\tau^{-1}I_{n-1},
\]
the system can be rewritten as
\[
\left(1+\Delta t\,\tau^{-1}\right)\Delta I_n
=
\Delta t\,b_n .
\]

For sufficiently small \(\Delta t\), namely if
\[
\Delta t\,\tau^{-1} \ll 1,
\]
the geometric series converges and the algebraic operator admits the Neumann expansion
\[
\left(1+\Delta t\,\tau^{-1}\right)^{-1}
=
1
-\Delta t\,\tau^{-1}
+\Delta t^2\,\tau^{-2}
-\cdots .
\]

Therefore,
\[
\Delta I_n
=
\left(
\Delta t
-\Delta t^2\,\tau^{-1}
+\Delta t^3\,\tau^{-2}
-\cdots
\right)b_n .
\]

This expression explicitly shows that the transient response is generated through repeated applications of the operator
\[
\tau^{-1}=L^{-1}R
\]
to the excitation term.

This Krylov-based interpretation naturally extends to the full-order MQS electromagnetic problem. After time discretization, the transient dynamics at each step are governed by a linear system whose dominant evolution modes are associated with the characteristic time constants of the generalized eigenvalue problem
\[
LI=\lambda RI.
\]

Equivalently, this relation can be written as
\[
R^{-1}LI=\lambda I,
\qquad
L^{-1}RI=\lambda^{-1} I.
\]

Since these two operators are the inverse of each other, they share the same eigenvectors and describe identical dynamical eigenspaces of the transient system.

From a computational viewpoint, however, the formulation based on \(R^{-1}L\) is significantly preferable. It avoids repeated inversions of the dense inductance matrix \(L\), which would otherwise become prohibitive for large-scale integral electromagnetic problems. This choice represents a key aspect of the proposed MOR strategy, since the reduced space is entirely constructed through successive applications of the operator \(R^{-1}L\), namely the evaluation of the direct integral operator followed by the inversion of \(R\) entails a negligible computational cost.

Finally, a direct circuit interpretation of this Krylov process is reported in Appendix~\ref{app:Krylov_circuit}.

While the Krylov approximation efficiently captures the dominant dynamical behavior of the system, the MOR framework also requires a a unified representation of the forcing sources. In practical fusion applications, the forcing term evolves during the transient because of both the time variation of the external sources and the recursive contribution of the previous state \(I_{n-1}\). Therefore, especially for rapidly varying plasma events, the MOR construction requires a low-dimensional representation of the forcing manifold capable of retaining the dominant transient dynamics without explicitly depending on the full time-marching evolution.
\section{Two-stage MOR construction}
\label{sec:two_stage_mor}

The MOR procedure is organized into two stages. First, the time-dependent excitations are compressed into a time-independent RHS matrix \(B_{\mathrm{mor}}\). Second, a reduced basis is generated from \(B_{\mathrm{mor}}\) through a source-driven block Krylov procedure.

An ideal reduced-order strategy would rely exclusively on the geometrical source operators, thus producing a fully static and event-independent reduced basis that could be reused for arbitrary transient scenarios. However, the numerical results presented in Section~\ref{sec:numerical_results} show that, for plasma transient events characterized by strongly varying temporal excitations, such a purely static construction leads to a rapid growth of the reduced-space dimension in order to preserve the desired accuracy. As a consequence, the resulting MOR model may become excessively large and computationally inefficient.

For this reason, temporal information must be incorporated into the forcing manifold through a suitable basis representation, allowing the MOR to retain only the dynamically relevant excitation directions. As shown by the numerical results in Section~\ref{sec:numerical_results}, this strategy significantly reduces the size of the reduced space with respect to a purely static construction while preserving the accuracy of the transient response.

\subsection{Stage 1: Construction of \(B_{\mathrm{mor}}\)}

The construction of \(B_{\mathrm{mor}}\) is performed independently for each excitation family through three sequential steps. First, the excitation waveforms are represented on a common multiresolution temporal basis composed of scaling functions and wavelets at different resolution levels. Second, a wavelet decomposition is applied to obtain a compressed coefficient matrix \(C\), retaining only the dominant transient contributions. Finally, the compressed forcing blocks are further reduced through a truncated singular value decomposition (SVD), yielding a orthonormal representation of the excitation.

The procedure can be summarized as
\begin{equation}
b(t)=V^{(k)}\,\alpha(t),
\qquad
\alpha_i(t)\approx
\sum_{m} a_{i,m}\,\phi(t-t_m)
+
\sum_{\ell=1}^{p_w}
C_{i,\ell}\,
\Psi_{\ell}(t),
\label{eq:wavelet_expansion}
\end{equation}
where \(\alpha_i(t)\) denotes the \(i\)-th component of the excitation vector \(\alpha(t)\). The functions \(\phi(t-t_m)\) are the shifted scaling functions associated with the coarse-scale representation, while \(\Psi_{\ell}(t)\) denotes the wavelet basis functions. The multi-index \(\ell=(j,k)\) collects both the resolution level \(j\) and the corresponding shift index \(k\), such that
\[
\Psi_{\ell}(t)
=
\Psi_{j,k}(t)
=
\Psi\!\left(2^{j}t-t_{j,k}\right).
\]

The coefficients \(a_{i,m}\) and \(C_{i,\ell}\) represent the corresponding scaling and wavelet amplitudes, respectively, while the retained wavelet coefficients are assembled into the compressed matrix \(C\).

\begin{equation}
B^{(V^{(k)})} = V^{(k)}\,C,
\qquad
B^{(V^{(k)})} \approx U_k\,\Sigma_k\,W_k^T,
\qquad
\widehat{B}^{(V^{(k)})} = U_k.
\label{eq:svd_compression}
\end{equation}

The effectiveness of this procedure relies on two main properties. First, wavelets provide a localized time-frequency representation that is particularly suitable for plasma transients characterized by short energetic events and multiple time scales. Second, since the electromagnetic problem is linear and time-invariant, neighboring wavelet slots --- which mainly correspond to translated versions of the same mother function --- generate dynamically similar responses. Therefore, the sources exhibit strong redundancy and can therefore be represented through a highly compressed set of spatial forcing directions.

Finally, the orthonormal compressed blocks associated with the four excitation families introduced in \eqref{eq:full2} and \eqref{eq:forcing} are assembled into the global forcing manifold
\[
B_{\mathrm{mor}}
=
\bigl[
\widehat{B}^{(\mathrm{axi})}
\;\;
\widehat{B}^{(\mathrm{3D})}
\;\;
\widehat{B}^{(\mathrm{volt})}
\;\;
\widehat{B}^{(J_0)}
\bigr].
\]

\subsection{Stage 2: Construction of the reduced basis}
\label{sec:mor_stage2}

The second stage of the MOR procedure consists in the construction of an excitation-driven reduced basis starting from the source manifold \(B_{\mathrm{mor}}\). Instead of approximating the full operator over the entire solution space, the reduced basis is generated only along the dynamical directions activated by the prescribed transient excitations.

Let \(A\) be the dynamic system matrix, defined as: \begin{equation}
A = L_K + \Delta t\,R_K .
\label{eq:matrix_A_def}
\end{equation}

When the null-space formulation is adopted, the projected operators \(L_K\) and \(R_K\) are defined according to~\eqref{eq:nullspace_continuous}.

\begin{figure}[H]
\centering
\resizebox{0.22\textwidth}{!}{%
\begin{tikzpicture}[
node distance=0.6cm,
>=latex,
block/.style={
rectangle,
draw,
rounded corners,
minimum width=4.2cm,
minimum height=0.65cm,
align=center
},
decision/.style={
diamond,
draw,
aspect=2,
align=center,
inner sep=2pt,
minimum width=2.5cm
},
line/.style={->, thick}
]


\node[block] (start) {
Initialize Krylov block\\
$W_0 = R^{-1}B_{\mathrm{MOR}}$
};

\node[block, below=of start] (basis0) {
Initialize basis\\
$V_0 = W_0$
};

\node[block, below=of basis0] (krylov) {
Generate Krylov block\\
$W_{k+1}=R^{-1}LW_k$
};

\node[block, below=of krylov] (proj) {
Projection step\\
$W_{k+1}\leftarrow W_{k+1}-V_k(V_k^TW_{k+1})$
};

\node[block, below=of proj] (update) {
Basis enrichment\\
$V_{k+1}=[V_k,W_{k+1}]$
};

\node[block, below=of update] (solve) {
Compute MOR residual\\
$\eta_{k+1}^{\mathrm{MOR}}$
};

\node[decision, below=0.6cm of solve] (check) {
$\eta_{k+1}^{\mathrm{MOR}} < \varepsilon_{\mathrm{MOR}}$
};

\node[block, below=0.6cm of check] (stop) {
Reduced basis\\
$V_{\mathrm{r}}$
};


\draw[line] (start) -- (basis0);
\draw[line] (basis0) -- (krylov);
\draw[line] (krylov) -- (proj);
\draw[line] (proj) -- (update);
\draw[line] (update) -- (solve);
\draw[line] (solve) -- (check);

\draw[line] (check) -- node[right, yshift=0.1cm]{yes} (stop);

\draw[line] (check.east) -- ++(1.5,0)
node[midway,above]{no}
|- (krylov.east);

\end{tikzpicture}%
}
\caption{Flowchart of the excitation-driven Krylov enrichment procedure used for the construction of the reduced basis \(V_{\mathrm{r}}\).}
\label{fig:krylov_flowchart}
\end{figure}
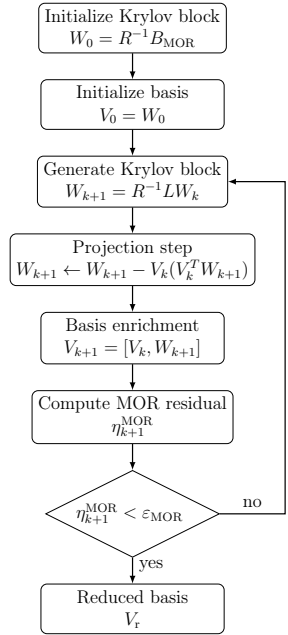

The reduced basis is then iteratively generated through a block Krylov enrichment procedure based on repeated applications of the operator \(R^{-1}L\) to the RHS. At each iteration, the newly generated Krylov block is orthogonalized with respect to the current reduced space and appended to the basis. The enrichment process continues until the prescribed MOR residual tolerance is satisfied. The overall workflow of the excitation-driven Krylov enrichment is summarized in Figure~\ref{fig:krylov_flowchart}.

The MOR tolerance \(\varepsilon_{\mathrm{mor}}\) is defined through the relative residual of the original full-order problem after projection onto the reduced space:
\begin{equation}
y_k = V_{\mathrm{r}} q_k,
\qquad
r_k = A V_{\mathrm{r}} q_k - B_{\mathrm{mor}},
\qquad
\eta_k^{\mathrm{mor}}
=
\frac{\|r_k\|_2}{\|B_{\mathrm{mor}}\|_2}.
\label{eq:mor_residual}
\end{equation}

\section{Reduced-Order Formulation and Computational Scalability}
\label{sec:reduced_order_models}

As previously discussed, the proposed MOR strategy does not make use of repeated explicit inversions or factorizations of the dense inductive operator \(L\) during basis generation. This property is confirmed by the Krylov construction procedure described in Section~\ref{sec:mor_stage2} and illustrated in Figure~\ref{fig:krylov_flowchart}. Indeed, the enrichment process only requires matrix--vector products of the form \(LI\), which can be efficiently evaluated through matrix-compression techniques such as hierarchical (\(H\)-matrix) methods, without explicitly assembling the full dense inductance matrix. The extension of the present framework toward compressed integral operators will be investigated in future developments.

When standard dense representations are adopted, the assembly of the inductive operator \(L\) remains computationally demanding due to the large-scale electromagnetic interactions involved, often requiring HPC resources for realistic fusion-device discretizations. This offline cost, however, is incurred only once and is effectively amortized over multiple transient simulations, since the same operator is reused for different plasma scenarios and operating conditions.

Moreover, because integral formulations discretize only conducting regions, localized geometrical modifications require recomputing only the corresponding interaction submatrices rather than rebuilding the entire global model. This property is particularly attractive for fusion-engineering applications, where uncertainties in geometry and material properties routinely motivate extensive parametric and optimization studies.

Once the reduced basis matrix
\[
V_{\mathrm{r}} \in \mathbb{R}^{N_{\mathrm{full}} \times N_{\mathrm{mor}}},
\qquad
N_{\mathrm{mor}} \ll N_{\mathrm{full}},
\]
has been constructed, where \(N_{\mathrm{mor}}\) denotes the dimension of the reduced space, the full-order null-space unknown is approximated as
\begin{equation}
y_n \approx V_{\mathrm{r}} q_n,
\label{eq:mor_ansatz}
\end{equation}
where \(q_n\in\mathbb{R}^{N_{\mathrm{mor}}}\) denotes the reduced state vector at time step \(n\).

Starting from the full-order null-space formulation
\begin{equation}
A y_n = b_n,
\label{eq:full_null_system}
\end{equation}
Galerkin projection onto the reduced space yields
\begin{equation}
V_{\mathrm{r}}^T A V_{\mathrm{r}} q_n
=
V_{\mathrm{r}}^T b_n,
\qquad
A_r = V_{\mathrm{r}}^T A V_{\mathrm{r}},
\qquad
b_{r,n} = V_{\mathrm{r}}^T b_n,
\label{eq:reduced_system}
\end{equation}
so that the reduced-order problem becomes
\begin{equation}
A_r q_n = b_{r,n}.
\label{eq:rom_system}
\end{equation}

After solving for \(q_n\), the full-order current vector is reconstructed as
\begin{equation}
I_n = K V_{\mathrm{r}} q_n + I_{0,n}.
\label{eq:rom_projection}
\end{equation}

The reduced system in~\eqref{eq:rom_system} preserves the algebraic structure of the full-order formulation introduced in~\eqref{eq:full_null_system}, with all operators and forcing terms replaced by their projected low-dimensional counterparts defined in~\eqref{eq:reduced_system}. The reduced approximation is introduced through~\eqref{eq:mor_ansatz}, while the reconstruction of the full-order currents follows from~\eqref{eq:rom_projection}.

The detailed algebraic derivation of the reduced operators is reported in Appendix~\ref{app:rom_definitions}, while the time discretization and numerical solution strategies directly follow the full-order methodology summarized in Appendix~\ref{app:transient}.

Although the same structural considerations apply to both the full-order and reduced formulations, the much smaller reduced dimension
\[
N_{\mathrm{mor}} \ll N_{\mathrm{full}}
\]
substantially changes the computational cost of the online stage. In particular, the reduced setup phase scales as
\[
\mathcal{O}(N_{\mathrm{mor}}^3),
\]
while the direct time-marching stage scales as
\[
\mathcal{O}(N_{\mathrm{mor}}^2),
\]
thereby avoiding the original
\[
\mathcal{O}(N_{\mathrm{full}}^3)
\qquad \text{and} \qquad
\mathcal{O}(N_{\mathrm{full}}^2)
\]
dependencies of the full-order formulation.

Consequently, the proposed MOR framework provides substantial computational speedups, with the reduced modal solver becoming particularly efficient for long transient simulations. The numerical experiments presented confirm these theoretical scaling properties.

\section{Model Order Reduction for the Null-Field Problem}
\label{sec:null_field_application}

The need for a null-field configuration in tokamaks originates from the physics of plasma breakdown~\cite{hender2007}. Plasma initiation occurs through a Townsend avalanche process driven by the toroidal electric field and becomes more effective in regions where the magnetic-field component orthogonal to the toroidal direction is sufficiently small. The minimum electric field required for breakdown is given by
\begin{equation}
E_{\min} = \frac{1.25 \times 10^{4} p}{\ln\left(510 p L\right)},
\qquad
L = 0.25\, a_{\mathrm{eff}} \frac{B_T}{B_\perp},
\label{eq:E_min_connection}
\end{equation}
where \(L\) denotes the connection length, \(B_T\) is the toroidal magnetic field, and \(B_\perp\) represents the magnetic field component orthogonal to the toroidal direction.

This relation shows that reducing $B_\perp$ increases the connection length $L$, thereby lowering the electric field required to initiate plasma breakdown. Since the available toroidal electric field is limited in practical devices, achieving a null-field configuration becomes a necessary requirement for reliable plasma start-up.

This physical requirement naturally leads to a control problem: the external coil currents must be adjusted so that the orthogonal magnetic field component approaches zero ($B_\perp \approx 0$) inside a prescribed spatial region during the breakdown phase.

\begin{figure}[htbp]
    \centering
    \includegraphics[width=0.45\textwidth]{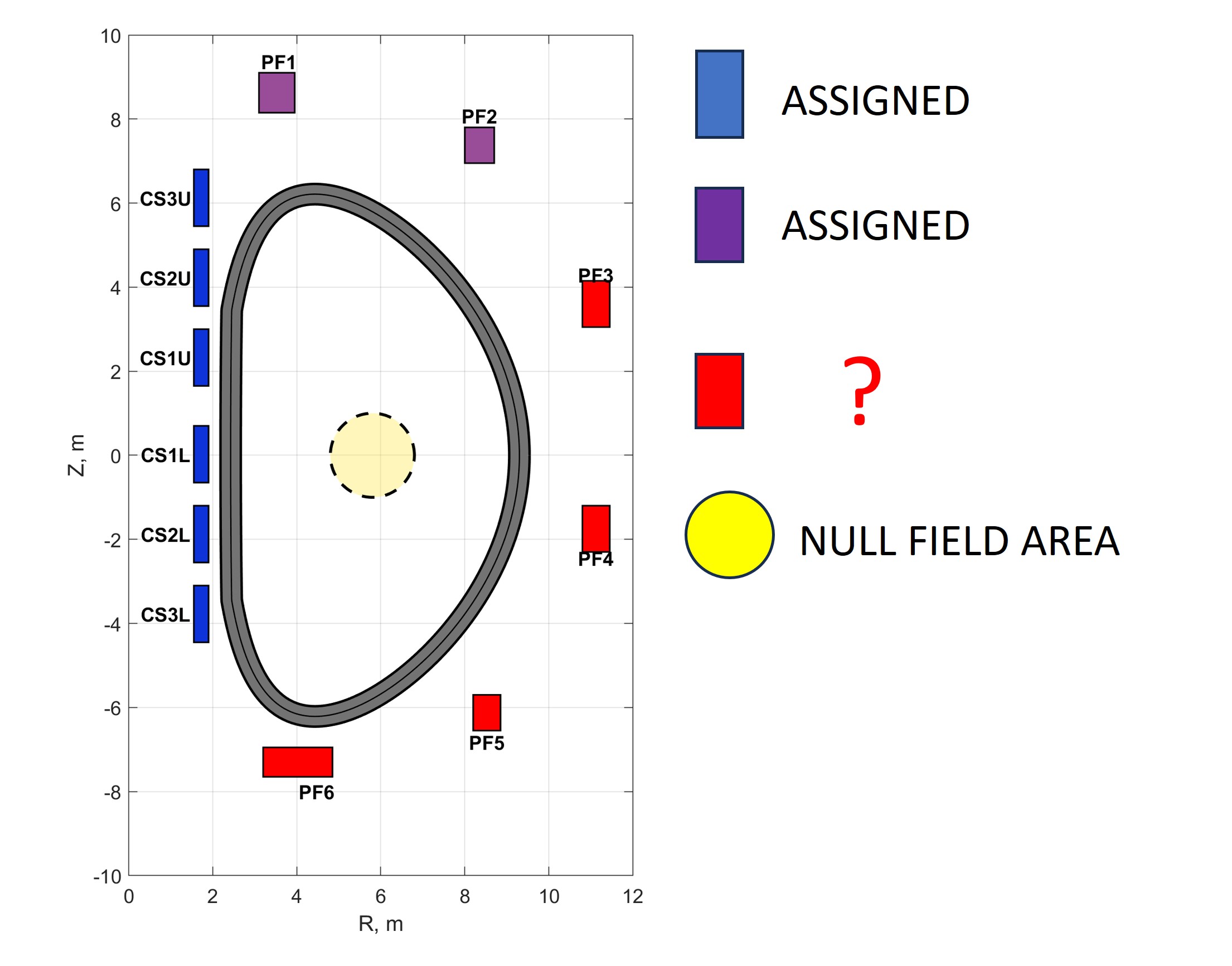}
    \caption{Simplified null-field problem: ITER poloidal section with CS and PF coils, highlighting the prescribed currents, the PF currents to be determined (red), and the region where the orthogonal magnetic-field component is required to vanish (yellow).}
    \label{fig:null_field_problem}
\end{figure}

The magnetic field in the region of interest results from the superposition of the contributions generated by the Central Solenoid (CS), the Poloidal Field (PF) coils, and the eddy currents induced within the passive conducting structures. In the considered configuration, the CS currents and the currents flowing in the first PF coils are prescribed, whereas the remaining PF-coil currents are treated as control variables and determined so as to enforce the null-field condition.

The present null-field configuration represents a simplified version of the complete plasma start-up control problem (see~\cite{digrazia2022waveforms} for additional details) and is adopted here as a benchmark configuration to assess the effectiveness of the proposed MOR methodology.

Accurate computation of the control currents requires a proper description of the transient electromagnetic response of the passive structures. While the present formulation is inspired by the approach proposed in~\cite{hl3_nullfield}, an important difference lies in the modeling of the conducting regions. In the cited work, only the vessel is represented through a simplified 2D axisymmetric discretization involving a limited number of conductors, whereas the present framework employs a detailed and fully resolved 3D representation of all passive structures.

Having introduced the null-field problem, the next two subsections describe the proposed acceleration strategy. First, the MOR technique speeds up the transient computations of the integral formulation. Second, this fast MOR framework generates the training dataset for a neural network surrogate, which then directly predicts the null-field solution practically in real-time.

\subsection{Integral formulation of the null-field problem}

Starting from the general full-order CARIDDI integral formulation introduced in~\eqref{eq:full2}, the null-field problem is obtained by retaining only the forcing contributions associated with the axisymmetric conductors. The null-space projection introduced in~\eqref{eq:nullspace_continuous} is then applied, leading to the coupled system

\begin{equation}
\begin{cases}
R_K \gamma(t)
+
L_K \dfrac{d\gamma}{dt}(t)
=
M_Y \dfrac{dY}{dt}(t)
+
M_X \dfrac{dX}{dt}(t),\\[10pt]
\mathbf{B}_{\perp}^{\mathrm{known}}(\gamma(t),Y(t))
+
\mathbf{D}_{\perp} X(t)
=
\mathbf{0},
\end{cases}
\label{eq:nullfield_coupled}
\end{equation}

where \(\gamma(t)\) denotes the null-space state vector associated with the induced currents, \(Y(t)\) collects the prescribed axisymmetric currents, including the CS coils and the first PF coils, and \(X(t)\) represents the unknown PF control currents. The term \(\mathbf{B}_{\perp}^{\mathrm{known}}\) contains the contribution to the orthogonal magnetic field generated by the prescribed sources and by the induced currents, while \(\mathbf{D}_{\perp}\) defines the mapping between the control currents and the orthogonal magnetic-field component within the null region.

The source term associated with the eddy-current dynamics is decomposed as

\[
\frac{d}{dt}\bigl(V_0 \alpha(t)\bigr)
=
V_{0Y}\frac{dY}{dt}(t)
+
V_{0X}\frac{dX}{dt}(t),
\]

with projected operators

\[
M_Y = K^T V_{0Y},
\qquad
M_X = K^T V_{0X}.
\]

The first equation in~\eqref{eq:nullfield_coupled} describes the transient electromagnetic dynamics projected onto the null-space basis, while the second equation imposes the cancellation of the orthogonal magnetic-field component within the target region.

The unknowns of the problem are therefore the electromagnetic state \(\gamma(t)\) and the control currents \(X(t)\), whereas the physical induced currents are reconstructed through

\[
I(t)=K\,\gamma(t).
\]

This formulation highlights that the computation of the PF control currents is intrinsically coupled with the transient electromagnetic response of the passive conducting structures. For completeness, the numerical solution procedures associated with both the full-order and reduced-order formulations are reported in Appendix~\ref{app:nullfield}.

\subsection{Neural Network Approach for Null-Field Control}
\label{sec:NN_null_field}

While MOR substantially reduces computational costs, repeated transient simulations may still remain too demanding for real-time digital-twin applications. To bypass this residual bottleneck, a data-driven Neural Network (NN) surrogate is introduced. The MOR solver is exploited offline as an efficient data-generation tool, while the NN learns the mapping between the excitation parameters and the resulting transient current dynamics.

As predicted by the Kolmogorov \(n\)-width theory~\cite{kolmogorov1956foundations}, constructing a surrogate model directly on raw, high-dimensional transient signals is impractical. To overcome this limitation, a POD--NN framework is designed to operate on a low-dimensional parameterization of both the inputs and outputs. First, restricting the excitation currents to typical ramp--plateau profiles allows the input waveforms to be uniquely described by a condensed parameter vector:
\begin{equation}
\boldsymbol{\mu}
=
\left[
A_{\mathrm{cs}},\,
t_{\mathrm{r,cs}},\,
A_{\mathrm{pf}},\,
t_{\mathrm{r,pf}},\,
t_{\mathrm{fin}}
\right]^T.
\label{eq:mu_param}
\end{equation}

Second, the generated control current responses exhibit a strong low-rank structure. They can therefore be compressed through a set of precomputed temporal POD modes \(\boldsymbol{\psi}_j(t)\), meaning the NN only needs to predict the reduced coefficients \(a_j^{\mathrm{NN}}(\boldsymbol{\mu})\). Once the offline training is completed, the online prediction of the transient response requires only a single forward NN evaluation followed by a fast linear combination:
\begin{equation}
x^{\mathrm{pred}}(t)
=
x_{\mathrm{mean}}(t)
+
\sum_{j=1}^{N_{\mathrm{POD}}}
a_j^{\mathrm{NN}}(\boldsymbol{\mu})
\,
\boldsymbol{\psi}_j(t).
\label{eq:NN_reconstruction}
\end{equation}

The comprehensive offline database generation, POD construction, network architecture, and training procedures are detailed in Appendix~\ref{app:POD_NN}.

%
\section{Numerical Results}
\label{sec:numerical_results}

The proposed MOR framework is assessed through a set of numerical benchmarks involving plasma events, plasma operating scenarios, and a simplified null-field problem. The analysis focuses on the approximation accuracy and computational efficiency achieved by the reduced-order formulation with respect to the corresponding full-order models. In addition, for the null-field application, the reduced-order model is exploited to construct and evaluate a neural-network surrogate model.

The reference configuration considered in this work is based on the ITER device and involves a large-scale unstructured discretization of the conducting structures relevant to electromagnetic transient analyses. The complexity of the model arises from the need to include realistic geometrical and physical details, whose electromagnetic response significantly affects the quantities of engineering interest.

\begin{figure}[H]
    \centering
    \includegraphics[width=0.92\textwidth]{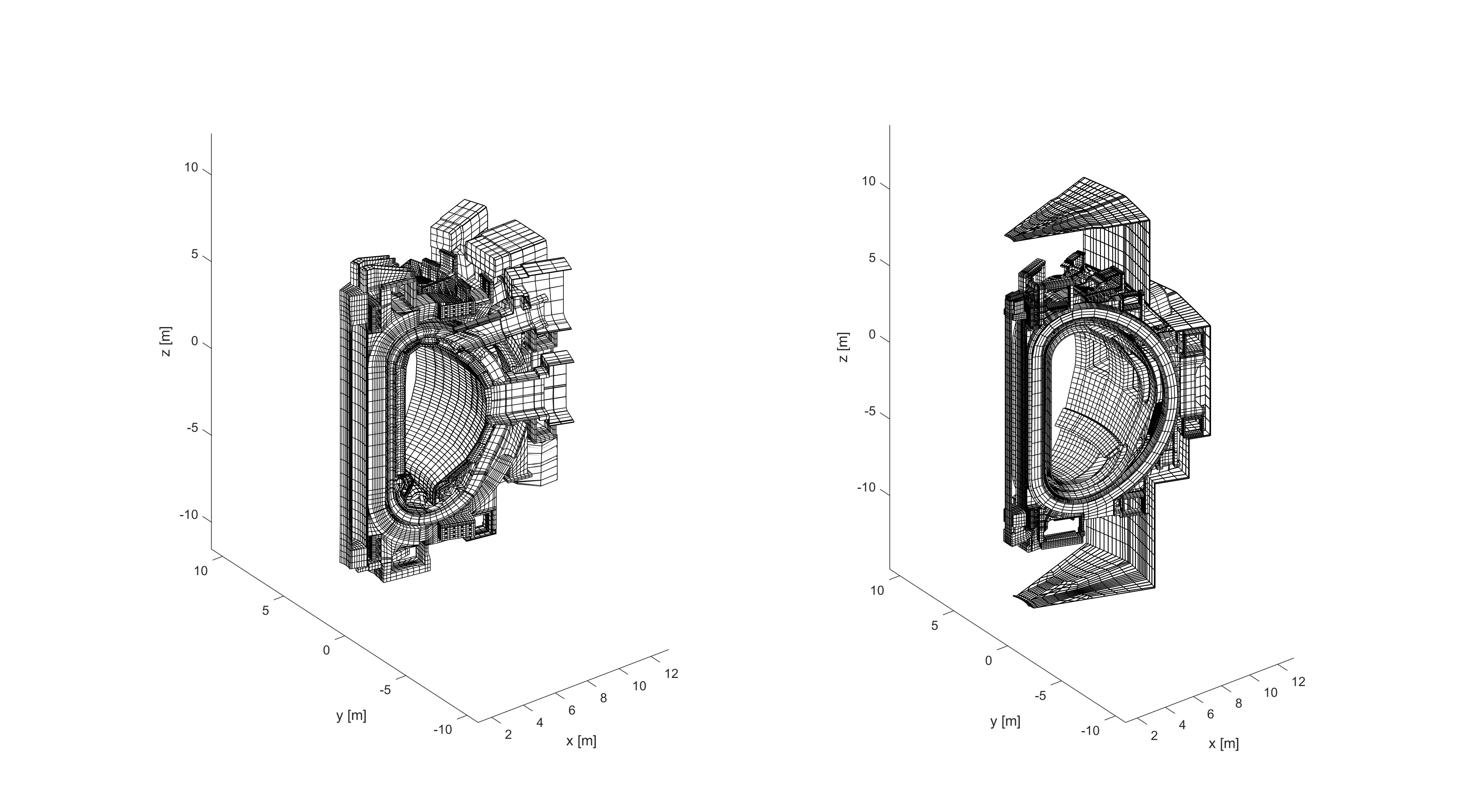}
\caption{40-degree sector model finite-element meshes adopted for the electromagnetic transient simulations: (left) Mesh\#1 used for plasma event analyses; (right) Mesh\#2 used for plasma scenarios and the null-field problem.}
\label{fig:two_meshes}
\end{figure}

Two computational meshes are considered in this work and are shown in Figure~\ref{fig:two_meshes}, while the corresponding discretization parameters and excitation configurations are summarized in Table~\ref{tab:mesh_data}. The excitation configurations are subdivided according to the four families of electromagnetic sources previously introduced in \eqref{eq:forcing} and \eqref{eq:full2}.

\begin{table}[H]
\centering
\caption{Finite Meshes and Excitations Information.}
\begin{tabular}{lcc}
\hline
 & Mesh \#1 & Mesh \#2 \\
\hline
Hexahedral elements & 80307 & 77410 \\
Tetrahedral elements & 0 & 3163 \\
Total elements & 80307 & 80573 \\
Mesh points & 147921 & 138762 \\
Degrees of freedom (DOFs) & 107742 & 106389 \\

Axis symmetric coils & 329 & 13  \\
3D coils & 80  & 0 \\
Voltage electrodes & 10  & 0 \\
Current electrodes & 7912 & 0 \\
\hline
\end{tabular}
\label{tab:mesh_data}
\end{table}

The two discretizations differ according to their intended applications. In particular, Mesh~\#1 is adopted for the analysis of electromagnetic loads generated during plasma transient events, with special emphasis on HALO currents, which are modeled through electrode representations. Mesh~\#2 is instead employed for AC-loss and subsequent thermal analyses under typical plasma-scenario conditions, where the cryostat is included in the model while electrode representations are neglected.

The numerical experiments were performed on three different computational platforms, whose main hardware characteristics are summarized in Table~\ref{tab:hardware_specs}. The first two systems are shared-memory workstations, whereas the third platform corresponds to the \emph{VEGA} partition of the ITER SDCC HPC cluster. The MOR, POD, and NN-training stages were executed on the same platforms in order to provide a consistent comparison of the computational performances.

\begin{table}[H]
\centering
\caption{Hardware specifications of the computational platforms used for the numerical experiments.}
\label{tab:hardware_specs}
\resizebox{0.98\textwidth}{!}{
\begin{tabular}{lccc}
\hline
\textbf{Specification} & \textbf{Machine \#1} & \textbf{Machine \#2} & \textbf{Machine \#3} \\
\hline
CPU model 
& Intel Xeon Gold 6338 @ 2.00\,GHz 
& AMD EPYC 7452 32-Core Processor 
& 2 $\times$ Intel Xeon Gold 6348 @ 2.60\,GHz \\

Physical CPU cores 
& 64 (2 $\times$ 32 cores) 
& 64 (2 $\times$ 32 cores) 
& 8 nodes $\times$ 56 cores/node \\

CPU frequency 
& 3.20\,GHz 
& 2.35\,GHz 
& 2.60\,GHz \\

RAM memory 
& 2.0\,TiB 
& 1.0\,TiB 
& 512\,GB DDR4 \\

CPU cache (L3) 
& 96\,MiB 
& 16\,MiB 
& 84\,MiB \\

GPU model 
& None
& 2 $\times$ NVIDIA Tesla V100S-PCIE-32GB 
& None \\

GPU memory 
& None
& 32\,GB each 
& None \\

Interconnect
& shared memory
& shared memory
& 2 $\times$ EDR InfiniBand, 100\,Gbit/s \\

Operating system 
& Ubuntu 22.04.4 LTS
& CentOS Linux 7 (Core)
& Red Hat Enterprise Linux 9.2 \\
\hline
\end{tabular}
}
\end{table}

In all numerical experiments, the reduced-order and reconstructed solutions are compared against the corresponding full-order models through suitable transient accuracy indicators and computational-performance metrics.

In particular, the accuracy of the MOR approximation is evaluated through the time-dependent relative error
\begin{equation}
\varepsilon_I(t)
=
\frac{\| I_{\mathrm{full}}(t)-I_{\mathrm{mor}}(t) \|}
{\| I_{\mathrm{full}}(t) \|}.
\label{eq:relative_error_current}
\end{equation}

This quantity measures, at each transient time instant, the discrepancy between the full-order and reduced-order solutions. Consequently, \(\varepsilon_I(t)\) represents a particularly stringent accuracy indicator, since it directly quantifies the instantaneous reconstruction error over the entire transient evolution.

\subsection{Plasma Events}

In this first numerical example, a set of four representative transient events typically occurring in tokamak devices is considered:
\[
\{\text{VDE\_UP}, \text{VDE\_DOWN}, \text{MD\_UP}, \text{MD\_DOWN}\}.
\] 

\begin{figure}[H]
    \centering
    \includegraphics[width=0.82\textwidth]{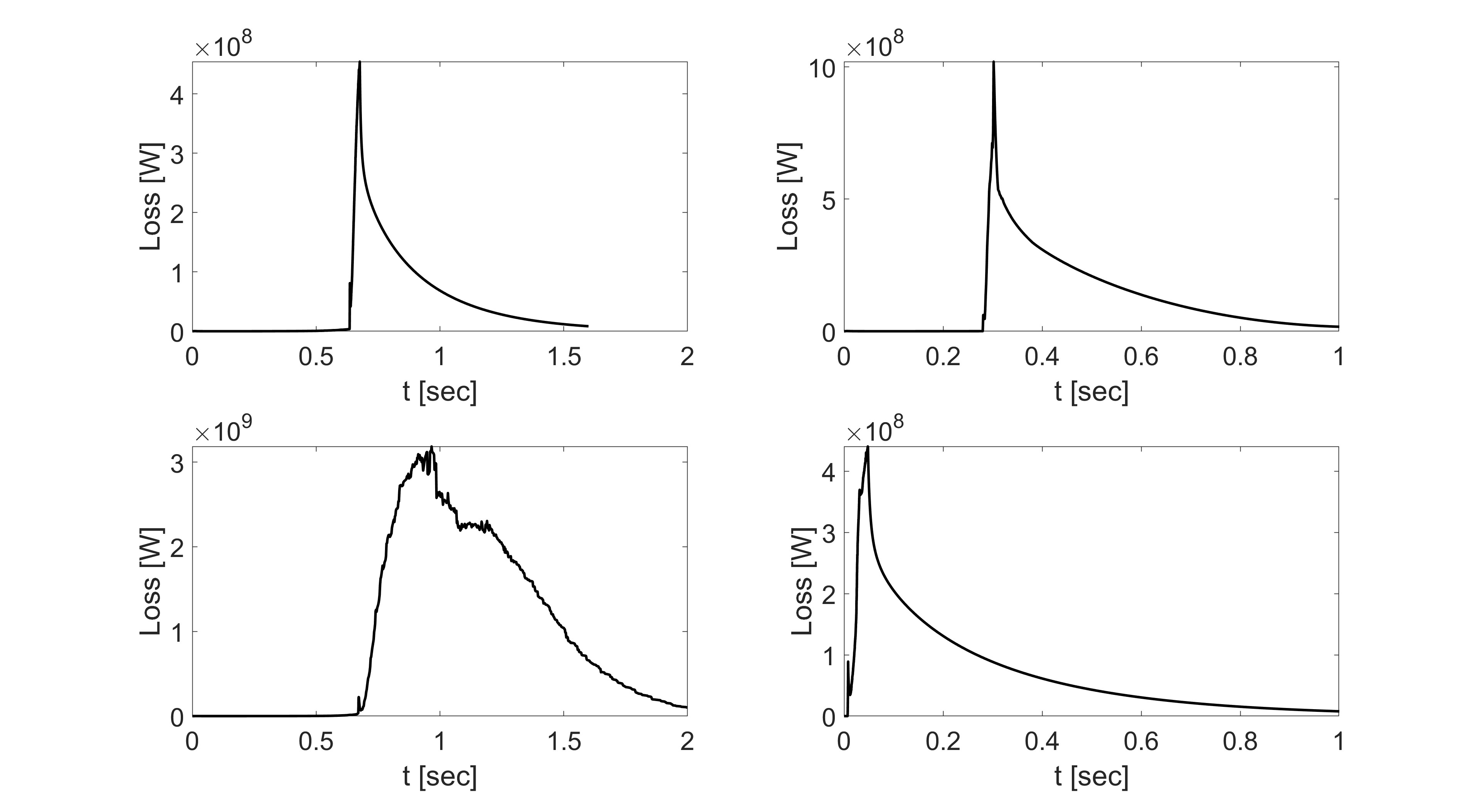}
    \caption{Power-loss evolution for the considered transient events: VDE\_UP (top-left), VDE\_DOWN (bottom-left), MD\_UP (top-right), and MD\_DOWN (bottom-right).}
    \label{fig:power_loss}
\end{figure}

Each of these events represents a specific, highly dynamic plasma instability. In particular, VDE\_UP and VDE\_DOWN (Vertical Displacement Events) are characterized by a rapid loss of vertical control, causing the plasma column to displace vertically toward the upper or lower regions of the vacuum vessel, respectively. On the other hand, MD\_UP and MD\_DOWN represent upward and downward Major Disruptions, which involve a sudden thermal crash followed by a fast current quench \cite{lehnen2015disruptions}. All these events induce intense electromagnetic transients, generating strong eddy currents and severe structural loads in the surrounding conducting components. As a consequence, they represent particularly demanding and highly relevant engineering benchmarks for fusion transient simulations. The temporal evolution of these phenomena is illustrated in Figure~\ref{fig:power_loss}, which reports the power-loss profiles associated with each considered transient event.

\begin{figure}[H]
    \centering
    \includegraphics[width=0.9\textwidth]{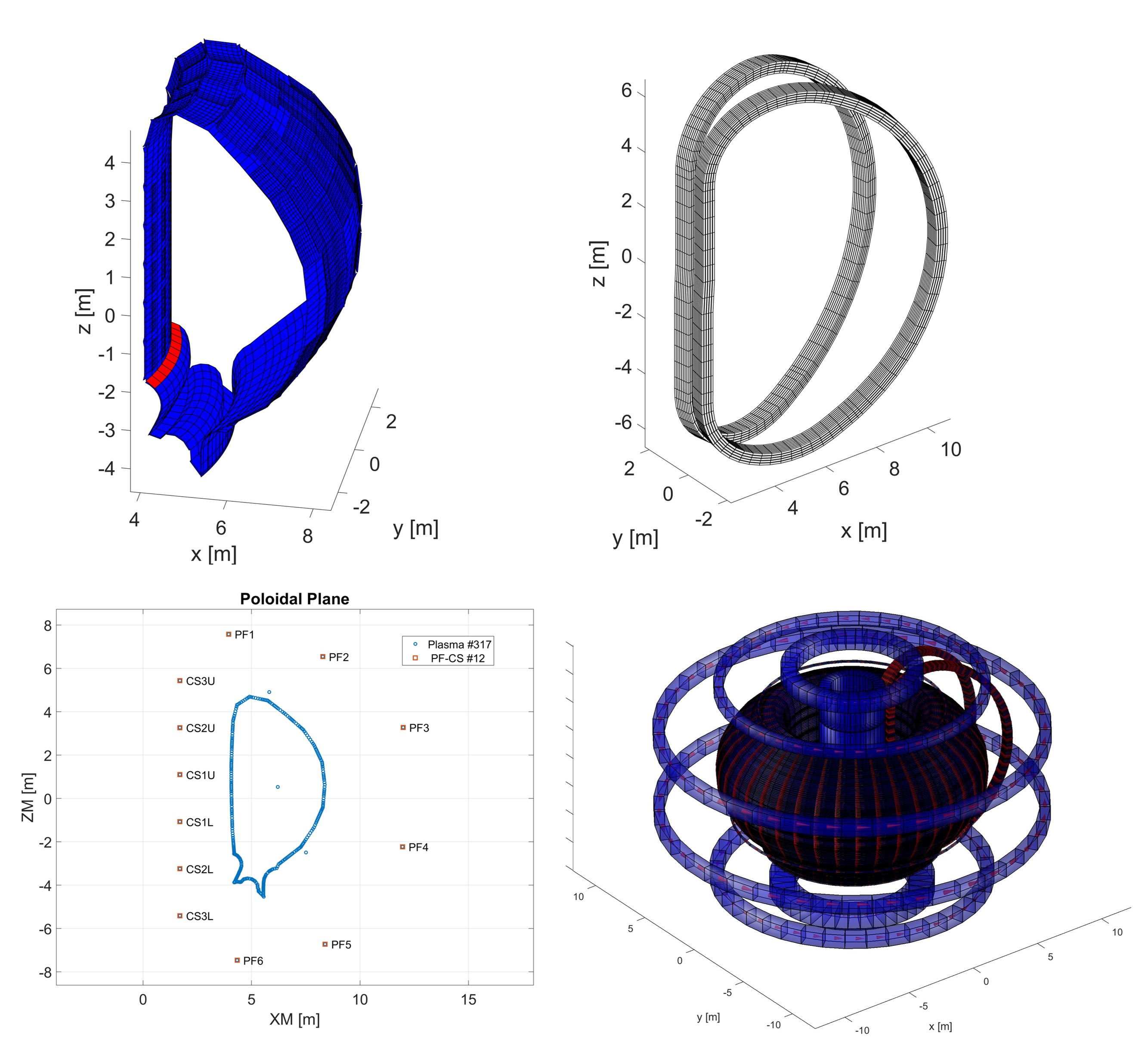}
    \caption{Top row: (left) voltage electrodes (red) and current electrodes (blue); (right) TF-coil mesh. Bottom row: (left) poloidal cross-section of the axisymmetric plasma, PF, and CS coil configuration; (right) 3D view of the complete active coil system.}
    \label{fig:axia_coil}
\end{figure}

The geometrical configuration of the excitation systems is depicted is depicted in Figure~\ref{fig:axia_coil}. This multi-family configuration introduces a strong spatial complexity. Furthermore, the excitation signals exhibit rapid temporal variations and highly irregular profiles, making these transients an exceptionally challenging benchmark for the MOR framework, which must accurately capture fast electromagnetic dynamics across a broad spectrum of time scales.

\begin{figure}[H]
    \centering
    \includegraphics[width=0.82\textwidth]{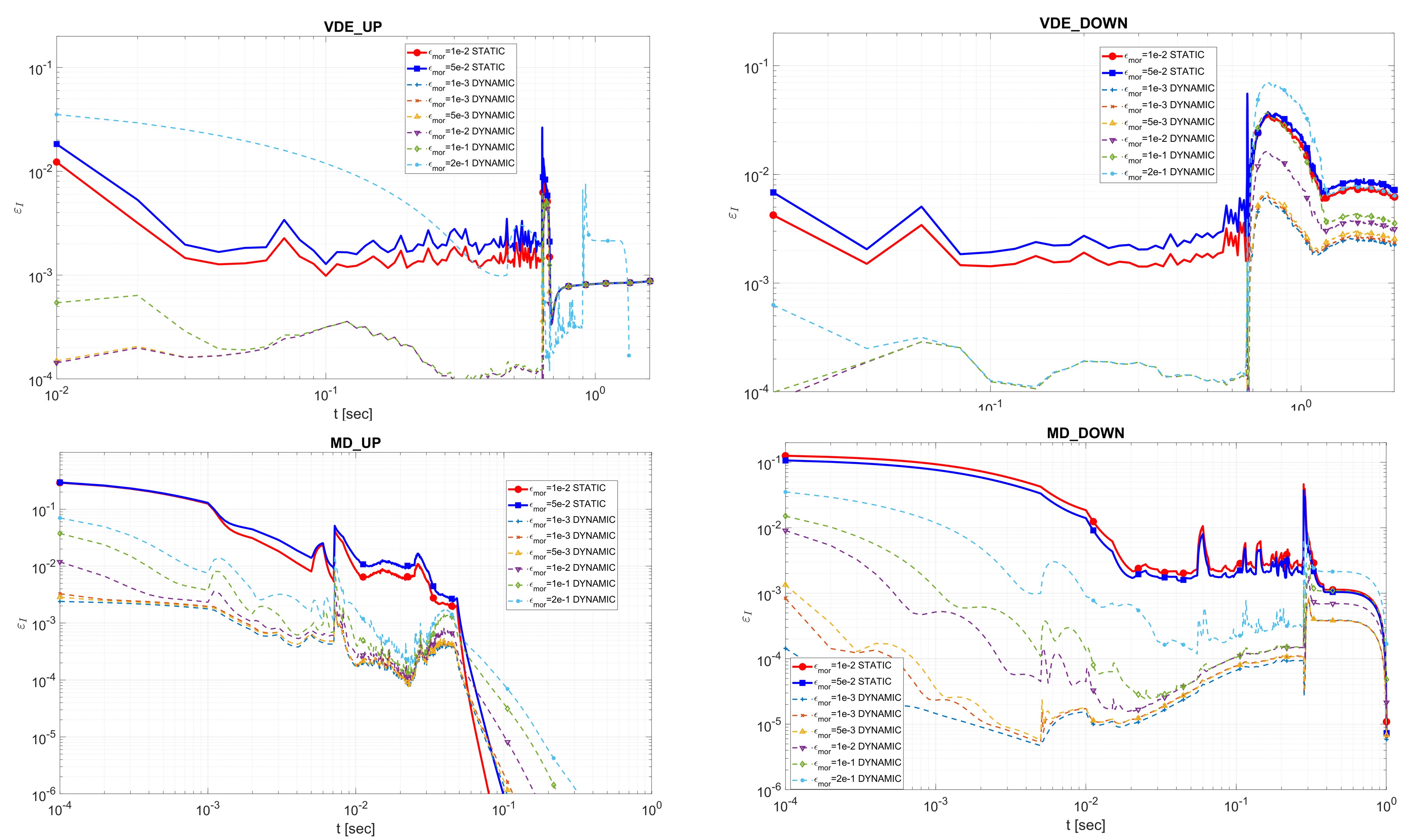}
    \caption{Time evolution of the relative error \(\varepsilon_I(t)\) for eight different values of \(\varepsilon_{\mathrm{MOR}}\) and four plasma-event configurations: VDE-UP (top-left), VDE-DOWN (top-right), MD-UP (bottom-left), and MD-DOWN (bottom-right).}
    \label{fig:error_di}
\end{figure}

\begin{figure}[H]
    \centering
    \includegraphics[width=0.82\textwidth]{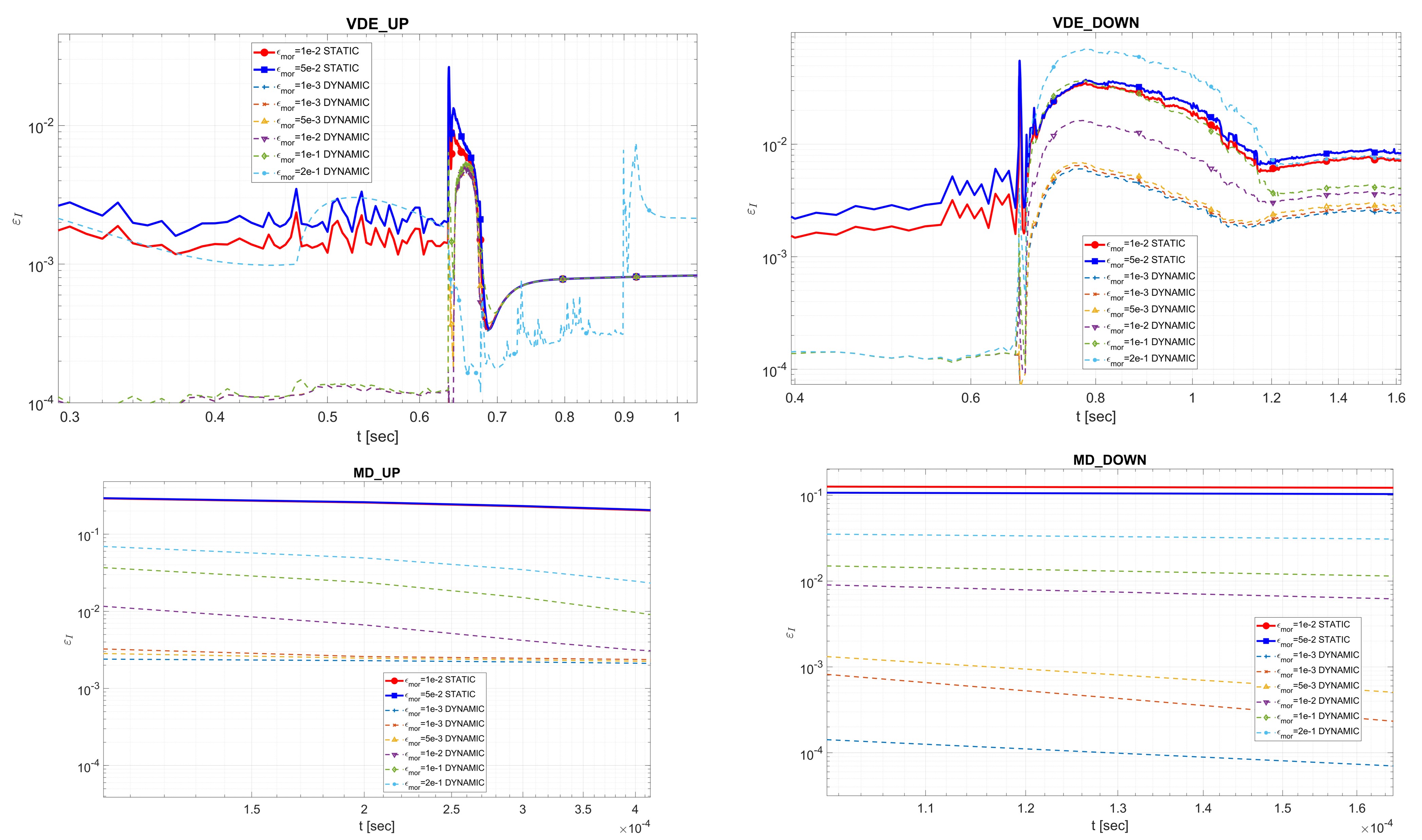}
    \caption{Zoomed view of the time evolution of the relative error \(\varepsilon_I(t)\) in the time interval where the maximum discrepancy occurs, for eight different values of \(\varepsilon_{\mathrm{MOR}}\) and four plasma-event configurations: VDE-UP (top-left), VDE-DOWN (top-right), MD-UP (bottom-left), and MD-DOWN (bottom-right).}
    
    \label{fig:error_di_zoom}
\end{figure}

Figure~\ref{fig:error_di} reports the time evolution of \(\varepsilon_I(t)\) for different truncation tolerances used in the compression of the forcing matrix \(B_{\mathrm{mor}}\). For clarity, an additional zoomed Figure ~\ref{fig:error_di_zoom} is also provided in order to highlight the time interval where the reconstruction accuracy is the lowest and the relative error reaches its maximum values.

To evaluate the efficiency of the proposed MOR strategy, the evolution of the reduced-order error is analyzed for different compression levels of the right-hand-side manifold employed in Stage~1 of the MOR construction.

\begin{table}[H]
\centering
\caption{MOR performance for the different plasma events.}
\label{tab:mor_performance}

\scriptsize
\resizebox{\textwidth}{!}{%
\begin{tabular}{m{1.8cm} c c c c c c c c c c}
\hline
Event &
$\varepsilon_{\mathrm{mor}}$ &
$N_{\mathrm{mor}}$ &
$N_{\mathrm{kry}}$ &
$\dfrac{N_{\mathrm{full}}}{N_{\mathrm{mor}}}$ &
$\left(\dfrac{N_{\mathrm{full}}}{N_{\mathrm{mor}}}\right)^2$ &
$\left(\dfrac{N_{\mathrm{full}}}{N_{\mathrm{mor}}}\right)^3$ &
Err. &
$T_{\mathrm{kry}}^{\mathrm{mor}}$ (s) &
$T_{\mathrm{trans}}^{\mathrm{mor}}$ (s) &
$\dfrac{T_{\mathrm{trans}}^{\mathrm{full}}}
{T_{\mathrm{trans}}^{\mathrm{mor}}}$ \\
\hline

\multirow{8}{=}{\centering\textbf{VDE\_UP}}
& $1\times10^{-2}$ & 19000 & 18 & 4.9 & $2.38\times10^{1}$ & $1.16\times10^{2}$ & $8.74\times10^{-3}$ & 1678 & 9922.71 & $1.5\times10^{1}$ \\
& $5\times10^{-2}$ & 13000 & 12 & 7.1 & $5.09\times10^{1}$ & $3.63\times10^{2}$ & $4.59\times10^{-2}$ & 1033 & 6211.30 & $2.4\times10^{1}$ \\
\cline{2-11}
& $1\times10^{-3}$ & 5781 & 40 & 16.0 & $2.57\times10^{2}$ & $4.13\times10^{3}$ & $9.53\times10^{-4}$ & 511 & 1361.66 & $1.09\times10^{2}$ \\
& $5\times10^{-3}$ & 4512 & 31 & 20.6 & $4.22\times10^{2}$ & $8.68\times10^{3}$ & $4.87\times10^{-3}$ & 388 & 749.73 & $1.99\times10^{2}$ \\
& $1\times10^{-2}$ & 3948 & 27 & 23.5 & $5.52\times10^{2}$ & $1.3\times10^{4}$ & $9.7\times10^{-3}$ & 363 & 575.49 & $2.59\times10^{2}$ \\
& $5\times10^{-2}$ & 2538 & 17 & 36.5 & $1.33\times10^{3}$ & $4.88\times10^{4}$ & $4.84\times10^{-2}$ & 233 & 304.38 & $4.89\times10^{2}$ \\
& $1\times10^{-1}$ & 1974 & 13 & 47.0 & $2.21\times10^{3}$ & $1.04\times10^{5}$ & $8.96\times10^{-2}$ & 226 & 190.73 & $7.81\times10^{2}$ \\
& $2\times10^{-1}$ & 1269 & 8 & 73.1 & $5.34\times10^{3}$ & $3.9\times10^{5}$ & $1.79\times10^{-1}$ & 159 & 83.45 & $1.79\times10^{3}$ \\
\hline

\multirow{8}{=}{\centering\textbf{VDE\_DOWN}}
& $1\times10^{-2}$ & 19000 & 18 & 4.9 & $2.38\times10^{1}$ & $1.16\times10^{2}$ & $8.74\times10^{-3}$ & 1021 & 10831.64 & $1.71\times10^{1}$ \\
& $5\times10^{-2}$ & 13000 & 12 & 7.1 & $5.09\times10^{1}$ & $3.63\times10^{2}$ & $4.59\times10^{-2}$ & 1004 & 7241.09 & $2.56\times10^{1}$ \\
\cline{2-11}
& $1\times10^{-3}$ & 8463 & 35 & 11.0 & $1.2\times10^{2}$ & $1.32\times10^{3}$ & $9.08\times10^{-4}$ & 1004 & 3512.64 & $5.28\times10^{1}$ \\
& $5\times10^{-3}$ & 6664 & 27 & 13.9 & $1.94\times10^{2}$ & $2.69\times10^{3}$ & $4.76\times10^{-3}$ & 829 & 2003.85 & $9.25\times10^{1}$ \\
& $1\times10^{-2}$ & 5950 & 24 & 15.6 & $2.43\times10^{2}$ & $3.78\times10^{3}$ & $8.84\times10^{-3}$ & 715 & 1599.47 & $1.16\times10^{2}$ \\
& $5\times10^{-2}$ & 3808 & 15 & 24.3 & $5.93\times10^{2}$ & $1.44\times10^{4}$ & $4.39\times10^{-2}$ & 466 & 703.67 & $2.63\times10^{2}$ \\
& $1\times10^{-1}$ & 2856 & 11 & 32.5 & $1.05\times10^{3}$ & $3.42\times10^{4}$ & $9.3\times10^{-2}$ & 362 & 445.71 & $4.16\times10^{2}$ \\
& $2\times10^{-1}$ & 1904 & 7 & 48.7 & $2.37\times10^{3}$ & $1.15\times10^{5}$ & $1.86\times10^{-1}$ & 293 & 228.13 & $8.12\times10^{2}$ \\
\hline

\multirow{8}{=}{\centering\textbf{MD\_UP}}
& $1\times10^{-2}$ & 19000 & 18 & 4.9 & $2.38\times10^{1}$ & $1.16\times10^{2}$ & $8.74\times10^{-3}$ & 1629 & 10179.06 & $1.46\times10^{1}$ \\
& $5\times10^{-2}$ & 13000 & 12 & 7.1 & $5.09\times10^{1}$ & $3.63\times10^{2}$ & $4.59\times10^{-2}$ & 1048 & 5806.16 & $2.57\times10^{1}$ \\
\cline{2-11}
& $1\times10^{-3}$ & 3400 & 49 & 27.3 & $7.44\times10^{2}$ & $2.03\times10^{4}$ & $9.46\times10^{-4}$ & 252 & 322.86 & $4.61\times10^{2}$ \\
& $5\times10^{-3}$ & 2652 & 38 & 35.0 & $1.22\times10^{3}$ & $4.27\times10^{4}$ & $4.84\times10^{-3}$ & 185 & 183.77 & $8.11\times10^{2}$ \\
& $1\times10^{-2}$ & 2380 & 34 & 39.0 & $1.52\times10^{3}$ & $5.91\times10^{4}$ & $8.65\times10^{-3}$ & 242 & 233.92 & $6.37\times10^{2}$ \\
& $5\times10^{-2}$ & 1564 & 22 & 59.3 & $3.51\times10^{3}$ & $2.08\times10^{5}$ & $4.56\times10^{-2}$ & 183 & 111.60 & $1.33\times10^{3}$ \\
& $1\times10^{-1}$ & 1156 & 16 & 80.2 & $6.43\times10^{3}$ & $5.16\times10^{5}$ & $9.87\times10^{-2}$ & 180 & 68.98 & $2.16\times10^{3}$ \\
& $2\times10^{-1}$ & 816 & 11 & 113.6 & $1.29\times10^{4}$ & $1.47\times10^{6}$ & $1.84\times10^{-1}$ & 134 & 34.68 & $4.3\times10^{3}$ \\
\hline

\multirow{8}{=}{\centering\textbf{MD\_DOWN}}
& $1\times10^{-2}$ & 19000 & 18 & 4.9 & $2.38\times10^{1}$ & $1.16\times10^{2}$ & $8.74\times10^{-3}$ & 1730 & 10310.86 & $1.44\times10^{1}$ \\
& $5\times10^{-2}$ & 13000 & 12 & 7.1 & $5.09\times10^{1}$ & $3.63\times10^{2}$ & $4.59\times10^{-2}$ & 990 & 5758.61 & $2.59\times10^{1}$ \\
\cline{2-11}
& $1\times10^{-3}$ & 5670 & 41 & 16.4 & $2.67\times10^{2}$ & $4.37\times10^{3}$ & $8.34\times10^{-4}$ & 520 & 1188.35 & $1.25\times10^{2}$ \\
& $5\times10^{-3}$ & 4320 & 31 & 21.5 & $4.61\times10^{2}$ & $9.89\times10^{3}$ & $4.99\times10^{-3}$ & 375 & 758.75 & $1.96\times10^{2}$ \\
& $1\times10^{-2}$ & 3915 & 28 & 23.7 & $5.61\times10^{2}$ & $1.33\times10^{4}$ & $8.45\times10^{-3}$ & 378 & 582.81 & $2.56\times10^{2}$ \\
& $5\times10^{-2}$ & 2430 & 17 & 38.2 & $1.46\times10^{3}$ & $5.56\times10^{4}$ & $4.81\times10^{-2}$ & 239 & 249.27 & $5.98\times10^{2}$ \\
& $1\times10^{-1}$ & 1890 & 13 & 49.1 & $2.41\times10^{3}$ & $1.18\times10^{5}$ & $9.08\times10^{-2}$ & 198 & 151.64 & $9.82\times10^{2}$ \\
& $2\times10^{-1}$ & 1215 & 8 & 76.3 & $5.82\times10^{3}$ & $4.44\times10^{5}$ & $1.89\times10^{-1}$ & 166 & 75.13 & $1.98\times10^{3}$ \\
\hline

\end{tabular}%
}
\end{table}

Table~\ref{tab:mor_performance} reports the performance of the proposed MOR
strategy for the four plasma events. The quantity $\varepsilon_{\mathrm{MOR}}$ denotes the prescribed stopping
tolerance adopted in the Krylov basis generation procedure described in the
flowchart of Figure~\ref{fig:krylov_flowchart}. The parameter $N_{\mathrm{kry}}$ represents the
number of Krylov iterations required to satisfy the prescribed tolerance,
while $N_{\mathrm{mor}}$ denotes the final dimension of the reduced subspace.

From columns 5 to 7, the reduction gains of the MOR formulation with respect
to the full-order problem are reported, highlighting the compression factors
obtained in terms of vector operations, dense matrix storage, and
cubic-complexity algebraic operations. The MOR accuracy is evaluated through the relative residual error
\(\eta_k^{\mathrm{mor}}\), introduced in \eqref{eq:mor_residual}.

Finally,
$T_{\mathrm{kry}}^{\mathrm{mor}}$
denotes the computational time required for the construction of the Krylov
reduced basis, while
$T_{\mathrm{trans}}^{\mathrm{mor}}$
represents the transient solution time of the reduced-order model. The last
column reports the transient speedup obtained by the MOR solver with respect
to the full-order transient simulation. All computational times reported in the following numerical this experiment are measured on Machine\#1. In particular, the corresponding full-order transient simulation requires approximately \(148980\,\mathrm{s}\).

The temporal discretization and wavelet-compression parameters adopted for the considered plasma events are summarized in Table~\ref{tab:event_time_parameters}. 
A uniform setup is employed for all cases, ensuring a consistent comparison of the MOR performance under different transient conditions.
\begin{table}[H]
\centering
\caption{Temporal discretization and wavelet parameters adopted for the considered plasma events.}
\label{tab:event_time_parameters}

\begin{tabular}{lccccc}
\hline
Event 
& \(t_{\mathrm{fin}}\) [s]
& \(\Delta t\) [s]
& \(\Delta t_{\mathrm{wavelet}}\) [s]
& Wavelet
& Max level
\\
\hline

VDE\_UP   
& 1.6 
& \(10^{-4}\) 
& \(10^{-3}\) 
& Daubechies-2 
& 1
\\

VDE\_DOWN 
& 1.99 
& \(10^{-4}\) 
& \(10^{-3}\) 
& Daubechies-2 
& 1
\\

MD\_DOWN  
& 1.6 
& \(10^{-4}\) 
& \(10^{-3}\) 
& Daubechies-2 
& 1
\\

MD\_UP    
& 1.6 
& \(10^{-4}\) 
& \(10^{-3}\) 
& Daubechies-2 
& 1
\\

\hline
\end{tabular}
\end{table}


The numerical experiments clearly demonstrate that, even under very strong compression of the forcing manifold, the reduced-order error remains remarkably low. This confirms the capability of the proposed MOR strategy to strongly compress the dynamical system dimension while accurately preserving the transient electromagnetic dynamics. As expected, smaller truncation tolerances enrich the forcing manifold and improve the approximation accuracy, at the cost of a moderate increase in the reduced-model dimension.

As previously anticipated, the numerical results confirm that static MOR strategies are intrinsically inefficient for rapidly varying transients. Despite using significantly larger reduced spaces, the static approach yields lower accuracy and higher costs than the proposed dynamic framework. For highly nonstationary electromagnetic transients, dynamically adapting the reduced basis to the excitation is therefore essential to simultaneously optimize precision and computational efficiency.

Finally, it is worth emphasizing that the measured transient speedup does not exactly coincide with the theoretical scaling factor \((N_{\mathrm{full}}/N_{\mathrm{mor}})^2\). The latter only represents the ideal asymptotic scaling associated with dense quadratic-complexity matrix operations, whereas the actual wall-clock time also includes implementation-dependent overheads, memory-access costs, data movement, and additional auxiliary computations that are not directly reduced by the dimensional compression.

\subsection{Plasma Scenario: verification of passive-structure effects on the ITER feeder busbar fields}

\begin{figure}[H]
    \centering
    \includegraphics[width=0.82\textwidth]{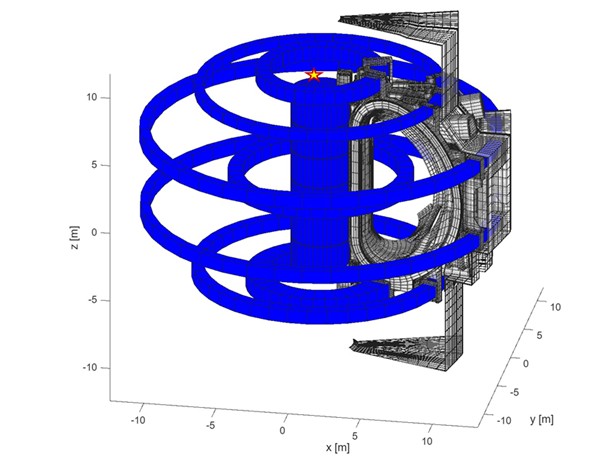}
  \caption{Mesh and excitation used in the ITER feeder busbar field computation.}
\label{fig:mesh_coils}
\end{figure}

To assess the impact of eddy currents induced in the passive structures (vacuum vessel, TF cases, etc.) on the local magnetic field within the ITER feeders, a comparative analysis was performed using the TRAPS~\cite{hertout2024traps} and CARIDDI codes, both validated for electromagnetic simulations of the ITER device. The TRAPS model neglects passive structures, whereas CARIDDI explicitly includes their electromagnetic coupling and induced eddy-current effects. The magnetic flux density was evaluated during the first \(100~\mathrm{s}\) of the ITER \(15~\mathrm{MA}\) plasma-pulse scenario at the terminal end of the CS1U busbar~1, whose location is marked with a star in Figure~\ref{fig:mesh_coils}, showing the computational mesh and the coil excitations.

In both simulations, the magnetic field generated directly by the busbars themselves was excluded from the calculation.

\begin{figure}[H]
    \centering
    \includegraphics[width=0.82\textwidth]{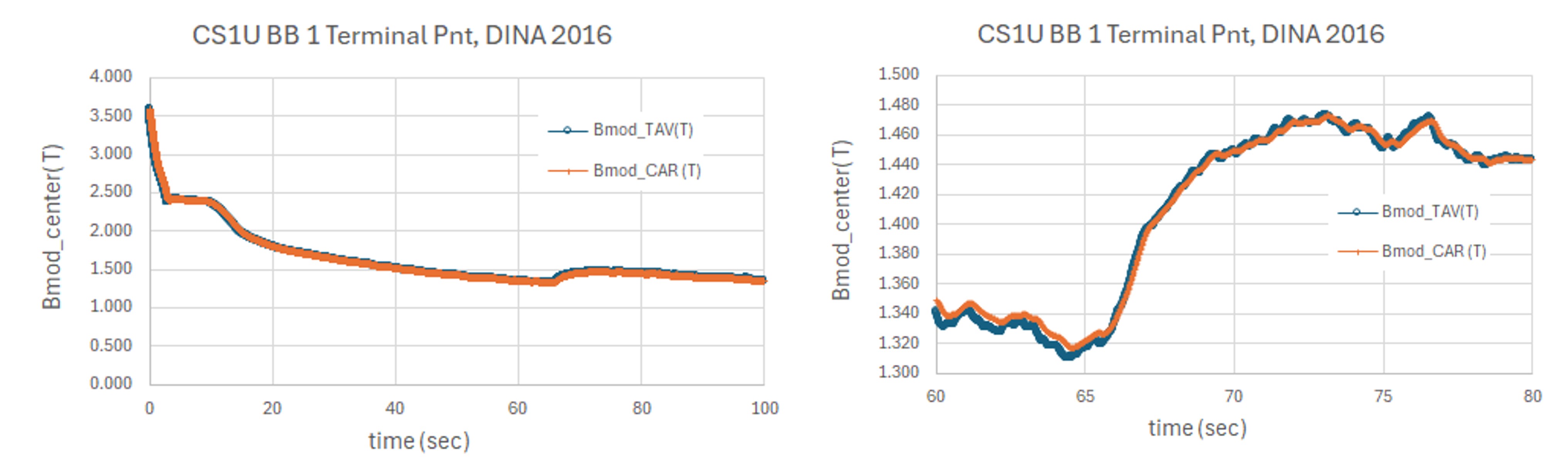}
  \caption{Flux-density distribution at the terminal end of the CS1U busbar center calculated with (left) TRAPS, neglecting passive structures and self-field effects, and (right) CARIDDI, including passive structures.}
\label{fig:cs1bb}
\end{figure}

The flux-density transient evaluated at the point
\[
(x=1.48006~\mathrm{m},\, y=-1.503~\mathrm{m},\, z=7.3225~\mathrm{m}),
\]
located at the terminal end of the CS1U busbar~1, is reported in Figure~\ref{fig:cs1bb}. A satisfactory agreement between the two models is observed. The CARIDDI solution clearly highlights the filtering effect introduced by the passive structures, whose induced eddy currents smooth the transient magnetic-field evolution. The contribution of the passive-structure eddy currents reaches a maximum value of approximately \(60~\mathrm{mT}\) during the initial discharge phase of the CS coils, corresponding to less than \(2\%\) of the total field.

\begin{figure}[H]
    \centering
    \includegraphics[width=0.82\textwidth]{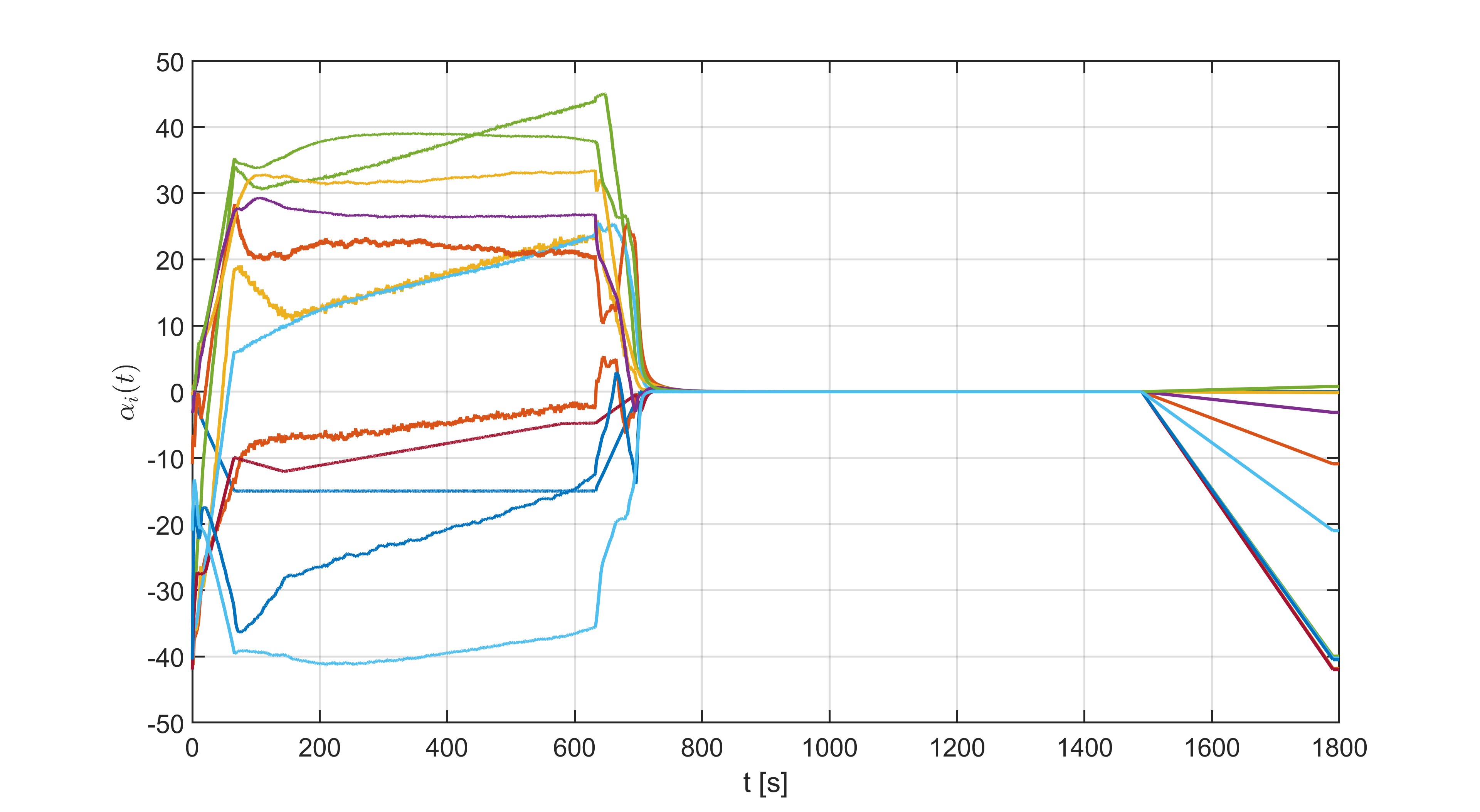}
  \caption{Time histories of the excitation currents associated with the 13 axisymmetric coils.}
\label{fig:scenario}
\end{figure}

The plasma-scenario simulations are carried out using Mesh\#2, where only the axisymmetric coil excitations are considered and no electrodes are included in the model. The corresponding time evolution of the excitation currents is reported in Figure~\ref{fig:scenario}.

Compared with plasma-event transients, the considered plasma scenario exhibits significantly smoother electromagnetic dynamics and more slowly varying excitation currents. Consequently, the MOR approximation becomes considerably more effective and requires only a limited reduced space to accurately reproduce the transient response. In this case, a static MOR construction is sufficient to achieve an accurate approximation while drastically reducing the computational cost.

The corresponding computational performances are summarized in Table~\ref{tab:traps_cariddi_mor}, where the proposed MOR approach is compared against the standard CARIDDI direct solver, the modal formulation, and the TRAPS code.

The reported results clearly demonstrate the remarkable efficiency of the MOR strategy. All computations reported in this comparison were performed on Machine\#3. While the standard CARIDDI transient simulation requires approximately \(96\) hours on \(64\) processors, and the modal approach about \(20\) hours, the MOR simulations are completed in only a few seconds on a single processor. Even for the most accurate reduced model, characterized by \(\varepsilon_{\mathrm{mor}}=10^{-2}\), the reduced dimension remains limited to only \(416\) unknowns.

The achieved accuracy is reported in Table~\ref{tab:eddy_errors}, where the relative errors on the magnetic flux-density components are shown for the considered MOR tolerances. In all cases, the discrepancy with respect to the eigenvalue solution remains very small, with errors typically ranging between \(10^{-5}\) and \(10^{-7}\).

\begin{figure}[H]
    \centering
    \includegraphics[width=0.82\textwidth]{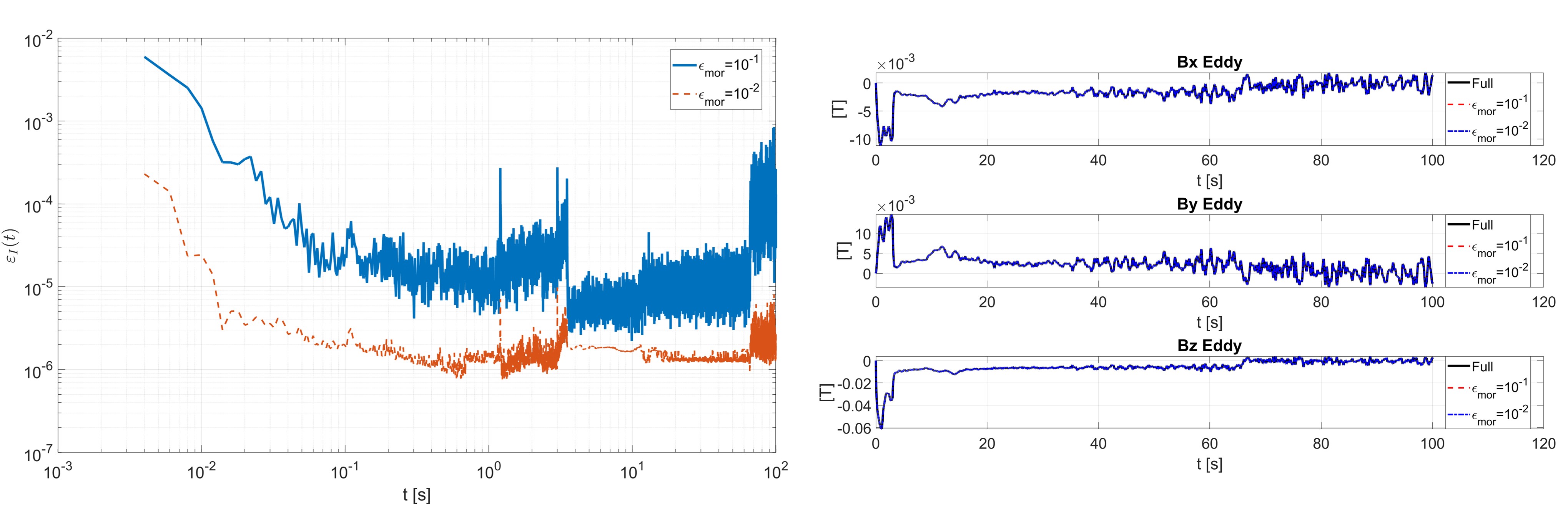}
\caption{Accuracy versus the MOR tolerance: (left) relative error on the transient solution, defined in~\eqref{eq:relative_error_current}; (right) relative error on the reconstructed magnetic flux density.}
\label{fig:scenario_accuracy}
\end{figure}

Finally, Figure~\ref{fig:scenario_accuracy} reports the influence of the MOR tolerance on both the transient solution accuracy and the reconstructed magnetic flux density. As expected, decreasing the MOR tolerance enriches the reduced basis and improves the approximation accuracy, while the computational cost remains negligible compared with the corresponding full-order simulations.

\begin{table}[htbp]
\centering
\caption{Comparison of computational times for the ITER feeder magnetic-field simulations.}
\label{tab:traps_cariddi_mor}

\begin{tabular}{lccc}
\hline
\textbf{Method} 
& \textbf{Processors} 
& \textbf{Dimension} 
& \textbf{Computational time}
\\
\hline

TRAPS 
& 1 
& -- 
& \(< 0.5\) h
\\

CARIDDI (Direct Method) 
& 64 
& 106389
& \( \sim 96 \) h
\\

CARIDDI (Modal) 
& 64 
& 106389 
& \( \sim 20 \) h
\\

CARIDDI--MOR 
\(\varepsilon_{\mathrm{mor}}=10^{-1}\)
& 1 
& 169 
& \( \sim 4 \) s
\\

CARIDDI--MOR 
\(\varepsilon_{\mathrm{mor}}=10^{-2}\)
& 1 
& 416 
& \( \sim 6 \) s
\\

\hline
\end{tabular}
\end{table}

\begin{table}[htbp]
\centering
\caption{Relative errors between the MOR and eigenvalue solutions for the eddy-current flux-density components.}
\label{tab:eddy_errors}

\begin{tabular}{lcc}
\hline
\textbf{Component}
& \(\boldsymbol{\varepsilon_{\mathrm{mor}}=10^{-1}}\)
& \(\boldsymbol{\varepsilon_{\mathrm{mor}}=10^{-2}}\)
\\
\hline

\(B_x\)
& \(1.6364\times10^{-5}\)
& \(6.0050\times10^{-7}\)
\\

\(B_y\)
& \(4.8268\times10^{-5}\)
& \(9.2320\times10^{-7}\)
\\

\(B_z\)
& \(1.3440\times10^{-5}\)
& \(2.8222\times10^{-7}\)
\\

\hline
\end{tabular}
\end{table}

\subsection{Null-field problem}
\label{sec:null_field_results}

This section presents the numerical validation obtained for the null-field problem. The objective is twofold: first, the capability of the proposed MOR strategy is assessed in terms of accuracy, computational speed, and compression capability; second, the reduced-order model is exploited as an efficient solver to generate the dataset required for the subsequent neural-network-based null-field framework. All computational timings reported in this numerical experiment are measured on Machine\#2.

\subsubsection{MOR performance for the null-field problem}
\label{subsec:null_field_mor}

\begin{figure}[htbp]
    \centering
    \includegraphics[width=0.85\textwidth]{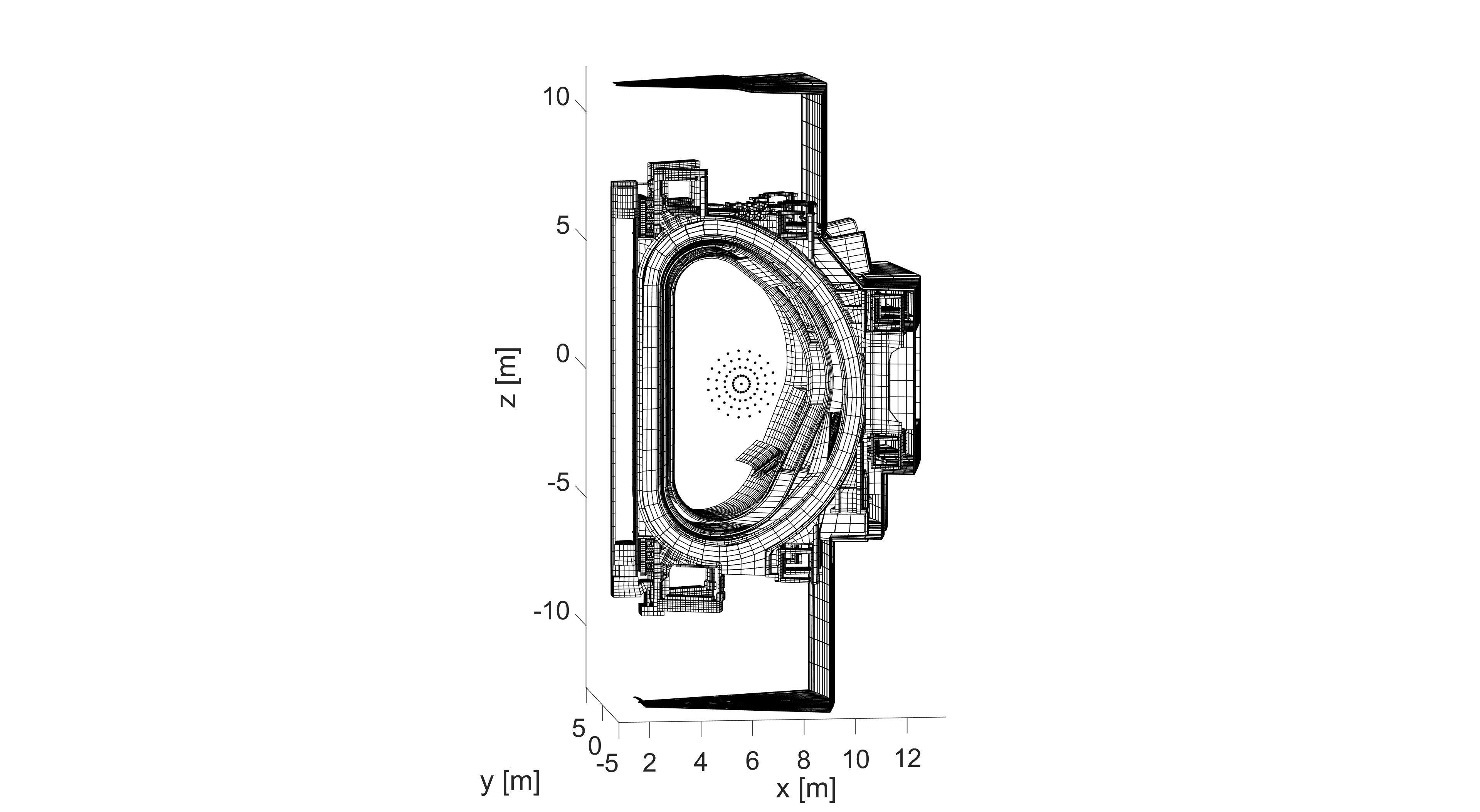}
    \caption{Computational mesh employed for the null-field problem, together with the location of the null-field control points.}
    \label{fig:null_field_mesh}
\end{figure}

The reduced-order model is initially validated by using the full-order solution as a reference. The MOR solution is obtained by projecting the governing equations onto the reduced space constructed through the excitation-driven Krylov procedure. These simulations are performed using Mesh\#2 (Table~\ref{tab:mesh_data}). As previously anticipated, a realistic full-3D discretization is adopted in order to accurately capture the electromagnetic coupling and the spatial distribution of the induced currents. Figure~\ref{fig:null_field_mesh} shows the corresponding computational mesh together with the control points employed for the poloidal-field cancellation.

Unlike the plasma-event —which are characterized by rapidly varying transients and therefore require larger reduced spaces to accurately reproduce the induced dynamics—the null-field configuration exhibits relatively smooth temporal evolutions. Consequently, a static MOR construction is sufficient in this case, leading to a significantly smaller reduced-order problem. Table~\ref{tab:mor_tolerance} summarizes the reduced-order metrics as a function of the prescribed MOR accuracy tolerance.

\begin{table}[H]
\centering
\caption{Influence of the MOR tolerance on reduced dimension, setup cost, and reconstruction accuracy.}
\begin{tabular}{cccccccc}
\hline
\(\varepsilon_{\mathrm{mor}}\) & Setup Time [s] & \(N_{\mathrm{it}}\) & \(N_{\mathrm{red}}\) 
& $\varepsilon_{\mathrm{PF3}}$
&
$\varepsilon_{\mathrm{PF4}}$
&
$\varepsilon_{\mathrm{PF5}}$
&
$\varepsilon_{\mathrm{PF6}}$
\\
\hline
\(10^{-1}\) & 128 & 2  & 24  & \(1.4\times10^{-3}\) & \(7.7\times10^{-3}\) & \(1.97\times10^{-2}\) & \(2.36\times10^{-2}\) \\
\(10^{-3}\) & 192 & 14 & 168 & \(2.0\times10^{-4}\) & \(1.6\times10^{-4}\) & \(2.0\times10^{-4}\)  & \(4.3\times10^{-4}\) \\
\(10^{-5}\) & 197 & 15 & 180 & \(2.0\times10^{-4}\) & \(1.6\times10^{-4}\) & \(2.0\times10^{-4}\)  & \(4.3\times10^{-4}\) \\
\hline
\end{tabular}
\label{tab:mor_tolerance}
\end{table}

For all null-field simulations,
\[
\Delta t = 10^{-3}\,\mathrm{s}.
\]

The accuracy of the null-field reconstruction is evaluated through the relative error on the reconstructed PF control currents \(X_i\), defined as
\begin{equation}
\varepsilon_{\mathrm{PF}i}
=
\frac{
\|
X_i
-
X_i^{\mathrm{mor}}
\|
}{
\|
X_i
\|
},
\qquad i=3,\dots,6.
\end{equation}

Table~\ref{tab:system_dimensions} compares the full-order and reduced-order models in terms of system dimension and memory occupation.

\begin{table}[htbp]
\centering
\caption{Comparison between the full-order and reduced-order models.}
\label{tab:system_dimensions}
\begin{tabular}{lccccc}
\hline
\textbf{Model} 
& \textbf{Dimension} 
& \textbf{Memory} 
& \(\mathbf{\dfrac{N_{\mathrm{full}}}{N_{\mathrm{mor}}}}\)
& \(\mathbf{\left(\dfrac{N_{\mathrm{full}}}{N_{\mathrm{mor}}}\right)^2}\)
& \(\mathbf{\left(\dfrac{N_{\mathrm{full}}}{N_{\mathrm{mor}}}\right)^3}\)
\\
\hline
FULL 
& 104,172 
& \(\sim 80\) GB 
& -- 
& -- 
& --
\\
MOR  
& 128--197 
& \(\sim 128\)--\(303\) kB 
& 529--814
& \(2.8\times10^5\)--\(6.6\times10^5\)
& \(1.5\times10^8\)--\(5.4\times10^8\)
\\
\hline
\end{tabular}
\end{table}

For the sake of clarity, the computation of the \(\mathbf{L}\) matrix, required by both formulations, was performed on the same machine using the FORTRAN-based CARIDDI code with 49 MPI processes, requiring approximately \(4022\)\,s.

\begin{figure}[H]
    \centering
    \includegraphics[width=0.85\textwidth]{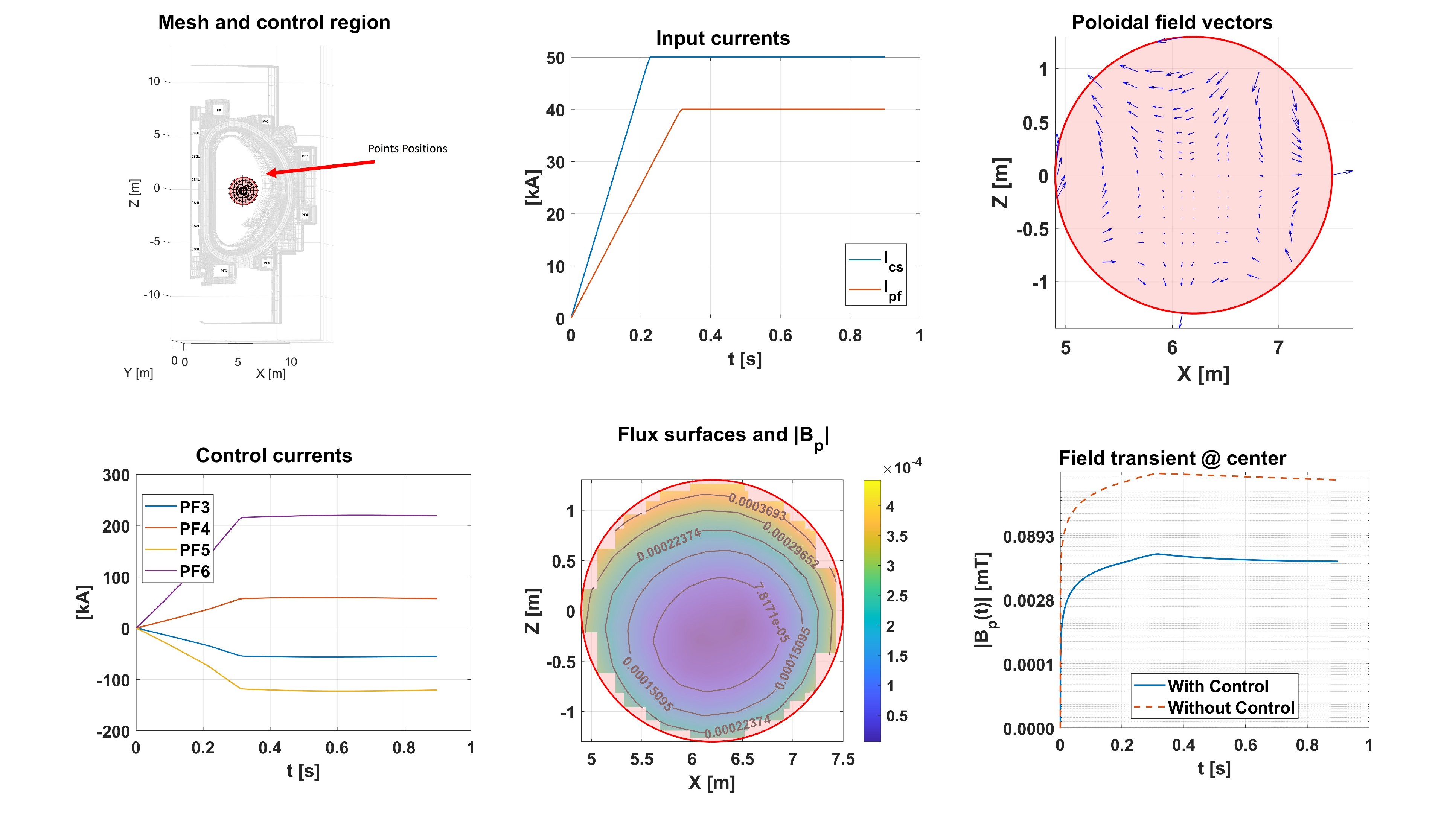}
\caption{
MOR transient performance for the ITER null-field problem.
(Top-left) ITER computational mesh and null-field region;
(Top-middle) CS and PF input currents;
(Top-right) poloidal magnetic field vectors;
(Bottom-left) obtained control PF currents;
(Bottom-middle) poloidal flux-density contour levels;
(Bottom-right) poloidal magnetic-field magnitude in the null-field region with and without active control during the transient.
}
  \label{fig:field_summary}
\end{figure}

Figure~\ref{fig:field_summary} provides a comprehensive overview of the transient simulation procedure for the ITER null-field problem. Specifically, the figure illustrates the computational mesh together with the null-field control region, the imposed CS and PF excitation currents, the resulting poloidal magnetic-field distribution, and the control PF currents computed through the reduced-order model. In addition, the corresponding magnetic flux-density contour levels and the temporal evolution of the poloidal magnetic-field magnitude at the center of the null-field region are reported for both the controlled and uncontrolled configurations.

The obtained control currents represent the dynamic values that must be applied in order to minimize the magnetic field at the prescribed control points during the transient evolution.

Overall, the obtained results demonstrate that the proposed MOR formulation achieves substantial reductions in computational cost and memory occupation while accurately preserving the transient null-field dynamics.

\subsubsection{NN Trained by MOR for the Null-Field Problem}
\label{subsec:NN_null_field_mor}

Table~\ref{tab:NN_summary} summarizes the main numerical settings and computational costs of the integrated MOR--POD--NN framework adopted for the null-field problem. The first block reports the nominal values and variation ranges used for the excitation parameters defining the transient input signals. The second block compares the computational costs of the FULL and MOR formulations, including setup times, solver times, MOR accuracy, and the average number of iterations required for convergence. Finally, the last block reports the training dataset size, the corresponding generation time, and the resulting neural network inference speed.

\begin{table}[H]
\centering
\caption{Summary of the MOR--POD--NN framework for the null-field problem.}
\label{tab:NN_summary}
\begin{tabular}{lc}
\hline
\textbf{Quantity} & \textbf{Value} \\
\hline
CS nominal amplitude \(A_{cs}^{nom}\) & \(45\,\mathrm{kA}\) \\
PF nominal amplitude \(A_{pf}^{nom}\) & \(50\,\mathrm{kA}\) \\
CS nominal rise time \(t_{r,cs}^{nom}\) & \(0.25\,\mathrm{s}\) \\
PF nominal rise time \(t_{r,pf}^{nom}\) & \(0.25\,\mathrm{s}\) \\
Nominal final time \(t_{\mathrm{fin}}^{nom}\) & \(1.0\,\mathrm{s}\) \\
Amplitude variation range & \([0.01,\,2]\times\) nominal value \\
Rise-time variation range & \([0.01,\,2]\times\) nominal value \\
Final-time variation range & \(0.7\text{--}1.0\,\mathrm{s}\) \\
\hline
MOR setup time & \(452.0\,\mathrm{s}\) \\
MOR relative error & \(3.386\times10^{-5}\) \\
MOR solver time & \(0.178\,\mathrm{s}\) \\
MOR average iterations per time step & \(14\) \\
FULL setup time & \(1556.8\,\mathrm{s}\) \\
FULL solver time & \(3486.0\,\mathrm{s}\) \\
\hline
Training database cardinality (\(N_{\mathrm{s}}\)) & \(3000\) \\
Dataset generation time & \(137.0\,\mathrm{s}\) \\
\hline
NN inference time & \(\approx 10^{-3}\,\mathrm{s}\) \\
\hline
\end{tabular}
\end{table}

\begin{figure}[H]
    \centering
    \includegraphics[width=0.78\textwidth]{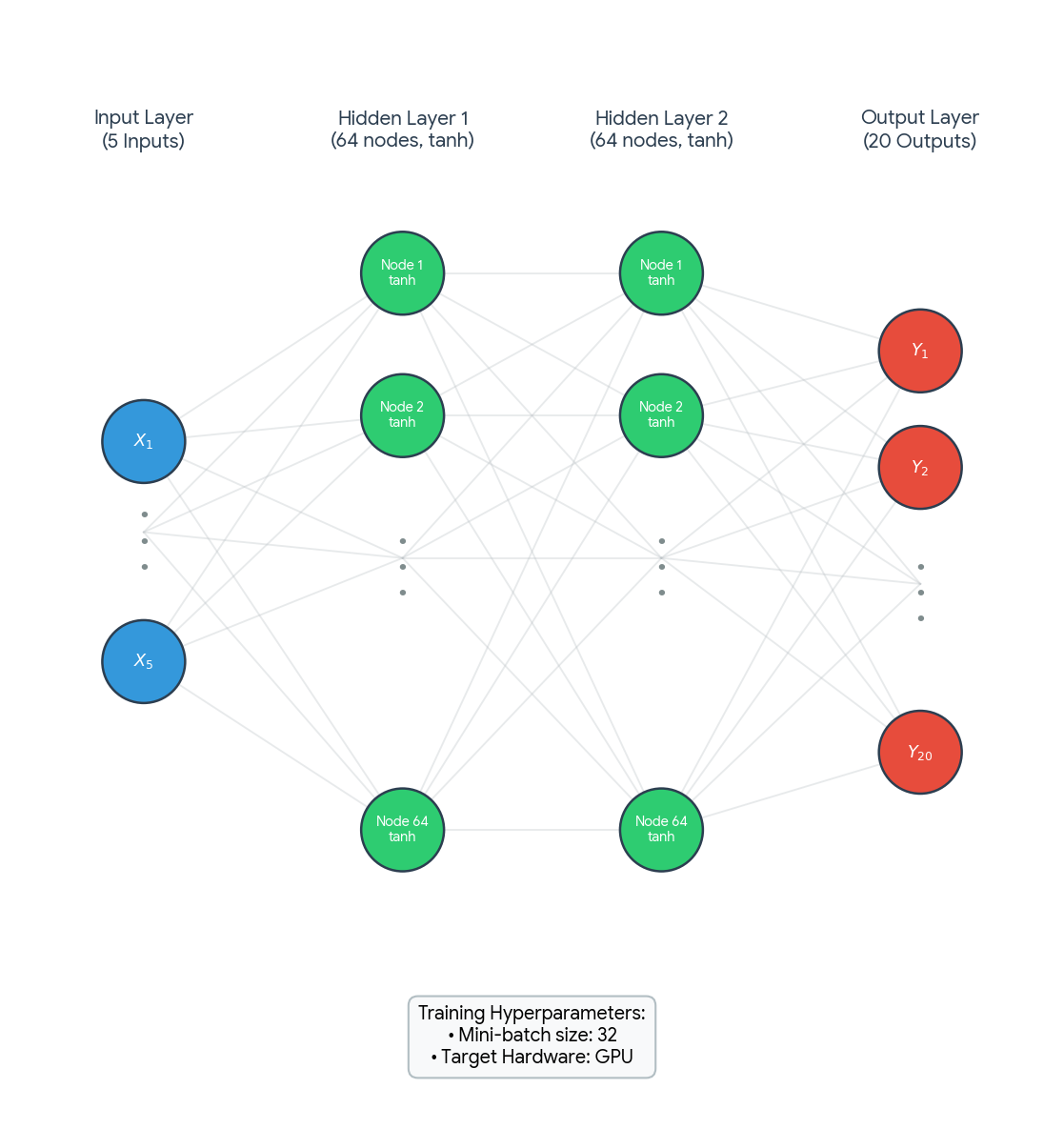}
    \caption{Architecture of the POD-NN surrogate model adopted in the numerical experiments. The network consists of three fully connected hidden layers with 64 neurons each and \(\tanh\) activation functions. The excitation parameter vector \(\mu\) is mapped onto the first 20 POD coefficients, which are then used to reconstruct the transient response through the POD basis functions.}
    \label{fig:nn_architecture}
\end{figure}

The POD-NN surrogate employed in the numerical experiments consists of a fully connected feed-forward neural network with three hidden layers, each composed of 64 neurons and hyperbolic tangent (\(\tanh\)) activation functions. 
The network receives the excitation parameter vector \(\mu\) as input and predicts the first 20 POD coefficients \(a_j(\mu)\), which are subsequently used to reconstruct the transient response through the POD basis functions. 
A schematic representation of the adopted architecture is reported in Figure~\ref{fig:nn_architecture}.

The results reported in Table~\ref{tab:NN_summary} highlight the remarkable computational efficiency of the proposed integrated approach. While the MOR transient solver requires less than \(0.2\,\mathrm{s}\) per simulation, the trained NN surrogate further accelerates the online phase by providing effectively instantaneous predictions of the control currents.

\begin{figure}[H]
    \centering
    \includegraphics[width=0.85\textwidth]{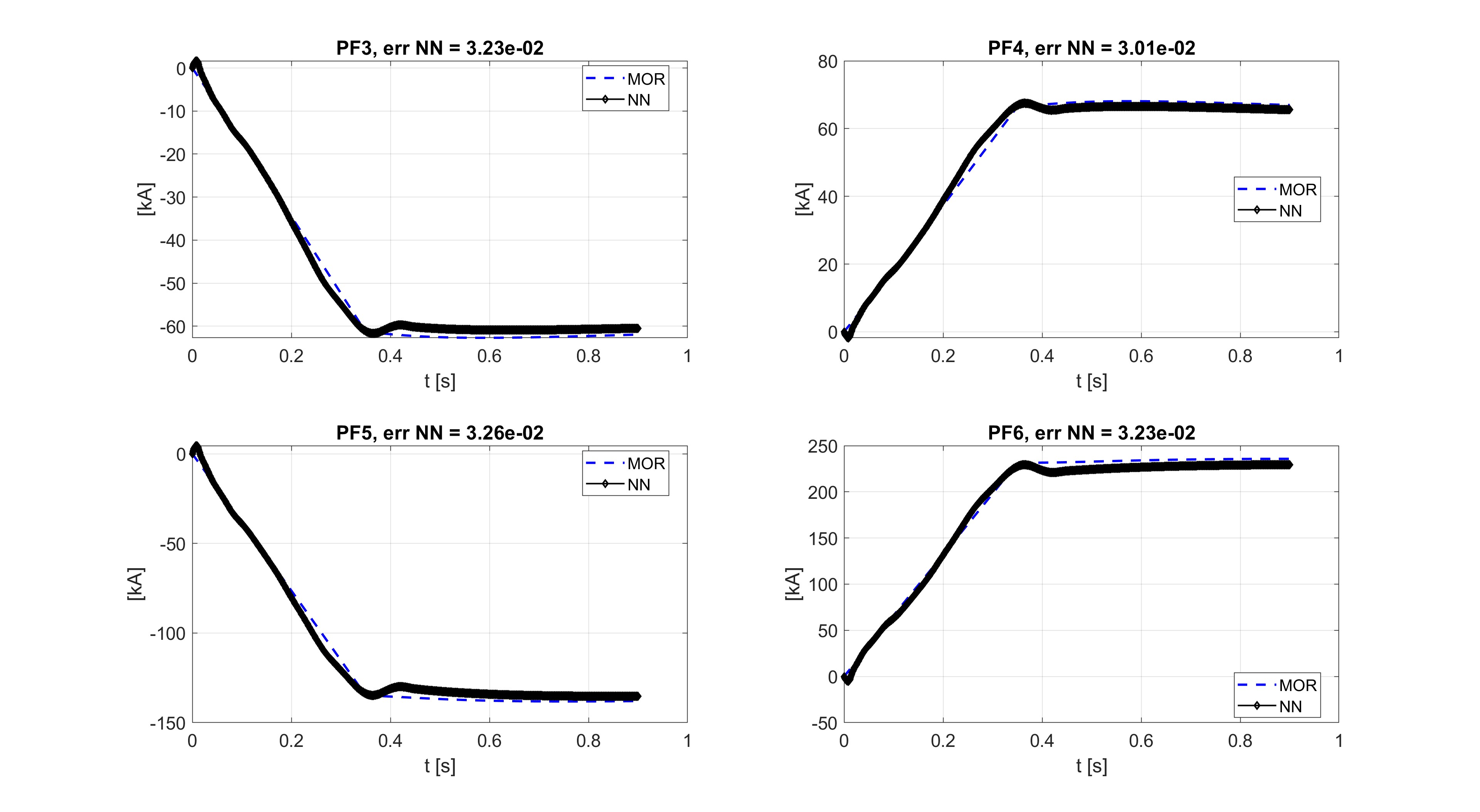}
    \caption{Comparison between the MOR and NN transient reconstructions of the control currents. For each PF current, the corresponding relative error of the NN prediction with respect to the MOR solution is also reported.}
    \label{fig:comparison_MOR_NN}
\end{figure}

Figure~\ref{fig:comparison_MOR_NN} compares the transient reconstructions obtained through the MOR solver and the trained neural-network surrogate. The NN reproduces the MOR transient behavior with good agreement across all PF currents.

By completely bypassing the transient numerical integration process, the trained NN surrogate enables sub-millisecond evaluations. These final metrics successfully demonstrate the suitability of the proposed MOR--POD--NN strategy for fast optimization loops and real-time inference in null-field control applications.

\section{Conclusion and Future Work}
\label{sec:conclusion}

This work presented an excitation-driven Model Order Reduction (MOR) strategy for large-scale electromagnetic transient problems arising from integral formulations in fusion-device applications. Unlike classical operator-driven approaches, the proposed methodology constructs the reduced space directly from the transient forcing manifold by combining wavelet-based temporal compression with source-driven Krylov projections. As a result, the reduced space is naturally tailored to the dynamically reachable electromagnetic responses associated with the considered excitation families.

The proposed framework was validated on several ITER-relevant plasma transient scenarios, including vertical displacement events, major disruptions, and a simplified null-field control problem. The numerical results demonstrated that the method achieves substantial reductions in system dimension while accurately preserving the underlying electromagnetic dynamics. The resulting computational speedups and memory savings make the approach particularly attractive for large-scale engineering analyses and repeated transient evaluations.

The application to the null-field problem further demonstrated how the MOR formulation can efficiently generate high-quality datasets for surrogate modeling. By combining MOR with POD and neural networks, real-time prediction of the control currents was achieved while maintaining high accuracy with respect to the full-order solutions. These results highlight the potential of the proposed methodology as a bridge between physics-based simulation tools and AI-enabled digital-twin architectures for fusion systems.

Future developments will focus on three main directions. First, the methodology will be extended toward fully self-consistent plasma--structure interaction models, where the electromagnetic response of passive structures is directly coupled with nonlinear plasma dynamics. Second, compressed representations of the inductive operator, including hierarchical matrices and matrix-free techniques, will be investigated in order to further improve scalability for extremely large-scale problems. Finally, more advanced AI-integration strategies will be explored, including adaptive reduced spaces, parametric MOR formulations, and physics-informed neural networks for real-time plasma control and optimization.

\appendix

\section{Cariddi integral formulation and discrete system}
\label{app:cariddi}

Under the magneto--quasi--static (MQS) approximation, the electromagnetic problem in a conducting non-magnetic domain \(V_c\) is described by expressing the electric field as~\cite{albanese1988,bossavit1981}
\[
\mathbf{E}
=
-\frac{\partial \mathbf{A}}{\partial t}
-
\nabla \phi,
\]
where \(\mathbf{A}\) is the magnetic vector potential and \(\phi\) is the scalar electric potential. The vector potential admits an integral representation in terms of the current density distribution within the conductor, namely
\[
\mathbf{A}(\mathbf{r},t)
=
\frac{\mu_0}{4\pi}
\int_{V_c}
\frac{\mathbf{J}(\mathbf{r}',t)}
{|\mathbf{r}-\mathbf{r}'|}
\, dV'
+
\mathbf{A}_s(\mathbf{r},t),
\]
where \(\mathbf{A}_s\) denotes the contribution generated by external sources. The constitutive relation is given by Ohm’s law
\[
\mathbf{J}=\sigma \mathbf{E}.
\]

Substituting the electric field into the constitutive relation and adopting a weak formulation leads to an integral equation for the current density, where the non-local nature of the magnetic interaction appears through the volume integral operator. The current density is automatically solenoidal by expanding it on the curl of a set of edge-based basis functions~\cite{bossavit1981,bossavit1998},
\[
\mathbf{J}(\mathbf{r},t)
=
\sum_{k=1}^{N}
I_k(t)\,
\nabla\times\mathbf{N}_k(\mathbf{r}),
\]
where \(I_k(t)\) are the unknown current coefficients. A Galerkin projection yields a system of differential-algebraic equations of the form
\[
L\frac{d\mathbf{I}}{dt}
+
R\mathbf{I}
+
F^T \boldsymbol{\Phi}
=
\frac{dV_0}{dt},
\]
where \(\boldsymbol{\Phi}\) collects voltages associated with current electrodes~\cite{rubinacci2002}.

The matrices entering the formulation are defined as follows. The resistive matrix is given by
\[
R_{ij}
=
\int_{V_c}
(\nabla\times\mathbf{N}_i)\cdot
\sigma^{-1}
(\nabla\times\mathbf{N}_j)
\, dV,
\]
and represents local dissipative effects (sparse, symmetric, and positive-definite matrix). The inductive matrix reads
\[
L_{ij}
=
\frac{\mu_0}{4\pi}
\int_{V_c}\int_{V_c}
\frac{
(\nabla\times\mathbf{N}_i)\cdot
(\nabla\times\mathbf{N}_j)
}
{|\mathbf{r}-\mathbf{r}'|}
\, dV\, dV',
\]
and accounts for long-range magnetic interactions (dense, symmetric, and positive-definite matrix).

Assuming that the source vector potential can be separated into spatial and temporal contributions as
\[
\mathbf{A}_s(\mathbf{r},t)
=
\mathbf{A}_0(\mathbf{r})\,I_s(t),
\]
its time derivative becomes
\[
\frac{\partial \mathbf{A}_s}{\partial t}
=
\mathbf{A}_0(\mathbf{r})
\frac{dI_s(t)}{dt}.
\]

For a single temporal waveform, the forcing contribution is written as
\[
V_{0,i}(t)
=
\int_{V_c}
(\nabla\times\mathbf{N}_i)\cdot
\frac{\partial\mathbf{A}_s}{\partial t}
\, dV
=
\frac{dI_s(t)}{dt}
\int_{V_c}
(\nabla\times\mathbf{N}_i)\cdot
\mathbf{A}_0(\mathbf{r})
\, dV.
\]

Defining the purely geometrical term
\[
V_{0,i}
=
\int_{V_c}
(\nabla\times\mathbf{N}_i)\cdot
\mathbf{A}_0(\mathbf{r})
\, dV,
\]
the forcing term assumes the factorized form
\[
V_{0,i}(t)
=
V_{0,i}\,\alpha(t),
\qquad
\alpha(t)
=
\frac{dI_s(t)}{dt}.
\]

Hence, the geometrical dependence is entirely contained in \(V_{0,i}\), whereas the temporal evolution is described by the scalar function \(\alpha(t)\). The extension to multiple independent temporal waveforms directly leads to the matrix formulation introduced in Eq.~\eqref{eq:forcing}.

The constraint matrix \(F\) enforces current conservation conditions on selected conductor surfaces and is defined as~\cite{rubinacci2002}
\[
F_{ij}
=
\int_{S_E}
(\nabla\times\mathbf{N}_j)\cdot
\mathbf{n}
\, dS,
\]
where \(S_E\) denotes the electrode surface and \(\mathbf{n}\) is the outward normal vector.

In the presence of geometric symmetries or periodicity conditions, additional algebraic constraints are introduced. For instance, symmetry planes or periodic sectors lead to relations of the form
\[
(F_{+\phi}-F_{-\phi})\,\mathbf{I}=0,
\]
where \(2\phi\) denotes the rotational symmetry angle of the considered cyclic sector. This condition ensures that the current distribution is consistent across matching boundaries.

\section{Full Model solution methods}
\label{app:transient}

The null-space formulation introduced in Section~\ref{sec:full_order_nullspace}, i.e., equation~\eqref{eq:nullspace_continuous}, is discretized in time through a backward Euler scheme. Let
\[
t_n=n\Delta t,
\qquad
y_n\approx y(t_n).
\]
The discrete problem reads
\begin{equation}
(L_K+\Delta t\,R_K)\,y_n
=
L_K\,y_{n-1}
+
\Delta t\,K^T b_{i,n}
-
\Delta t\,K^T R\,I_{0,n}
-
K^T L\,(I_{0,n}-I_{0,n-1}),
\label{eq:discrete_nullspace}
\end{equation}

Defining
\[
Z=L_K+\Delta t\,R_K,
\]
the system is solved at each time step as
\begin{equation}
Zy_n=b_n,
\label{eq:linear_system}
\end{equation}
with
\begin{equation}
b_n=
L_K\,y_{n-1}
+
\Delta t\,K^T b_{i,n}
-
\Delta t\,K^T R\,I_{0,n}
-
K^T L\,(I_{0,n}-I_{0,n-1}).
\label{eq:bn_term}
\end{equation}

Two alternative solution strategies are considered: a direct solver based on matrix factorization and a modal solver based on a generalized eigenvalue decomposition.

\paragraph{Direct solver.}
The matrix \(Z\) is factorized once through a Cholesky decomposition, with a computational cost proportional to \(\mathcal{O}(N_{\mathrm{full}}^3)\). Each subsequent time step then only requires forward and backward substitutions, leading to an online cost of order \(\mathcal{O}(N_{\mathrm{full}}^2)\).

\paragraph{Modal solver.}
Alternatively, the system is solved in a modal basis obtained from the generalized eigenvalue problem
\[
L_K\phi_k=\lambda_k R_K\phi_k,
\]
with the simultaneous diagonalization properties
\[
\Phi^T R_K\Phi=I_{\lambda},
\qquad
\Phi^T L_K\Phi=\Lambda,
\]
where
\(
\Lambda \in \mathbb{R}^{N_{\mathrm{full}}\times N_{\mathrm{full}}}
\)
is the diagonal matrix containing the generalized eigenvalues \(\lambda_k\), and
\(
I_{\lambda}
\in
\mathbb{R}^{N_{\mathrm{full}}\times N_{\mathrm{full}}}
\)
denotes the identity matrix.

Expanding the discrete solution vector in the modal basis as
\[
y_n=\Phi a_n,
\]
equation~\eqref{eq:discrete_nullspace} diagonalizes into
\[
(\Lambda+\Delta t\,I_{\lambda})a_n
=
\Phi^T b_n.
\]

The modal coefficients are therefore decoupled and can be updated independently at each time step as
\begin{equation}
a_{n,k}
=
\frac{(\Phi^T b_n)_k}{\lambda_k+\Delta t},
\qquad
k=1,\dots,N_{\mathrm{full}},
\end{equation}
and the full state is reconstructed through
\[
y_n=\Phi a_n.
\]

Since all modal equations are uncoupled, the transient integration becomes essentially linear in the modal dimension.

\section{Reduced Definitions and Solutions}
\label{app:rom_definitions}

The reduced-order formulation preserves the same algebraic structure of the full-order null-space problem. Starting from the reduced basis
\(
V_r \in \mathbb{R}^{N_{\mathrm{full}}\times N_{\mathrm{mor}}},
\)
the reduced operators are obtained through Galerkin projection:
\[
L_r
=
V_r^T L_K V_r,
\qquad
R_r
=
V_r^T R_K V_r,
\]
with
\[
L_K = K^T L K,
\qquad
R_K = K^T R K.
\]

Consequently, the reduced system matrix associated with the discrete transient dynamics becomes
\begin{equation}
\label{eq:reduced_system_matrix}
Z_r
=
L_r
+
\Delta t\,R_r.
\end{equation}

Similarly, the corresponding right-hand side forcing terms are projected onto the reduced space through left multiplication by \(V_r^T\). All operators and source terms are thus replaced by their lower-dimensional, projected counterparts. As a consequence, both the direct time-marching strategy and the modal formulation derived for the full-order problem remain formally valid and computationally accelerated within the reduced-order framework.







\section{Circuit Interpretation of the Krylov Expansion}
\label{app:Krylov_circuit}

An MQS electromagnetic system can be interpreted as the continuous counterpart of an RL network, where \(L\) and \(R\) denote the mutual inductance and resistance matrices, respectively.

Starting from an initial current distribution \(i_0\), the corresponding magnetic flux linkages are
\begin{equation}
\Phi_0 = L i_0 .
\end{equation}

In a resistive--inductive transient, the induced electromotive forces are generated by the temporal variation of the magnetic flux linkages. The resistive response of the network maps these induced voltages back into the current space through the action of the operator \(R^{-1}\), yielding the current redistribution
\begin{equation}
i_1 = R^{-1} L i_0 .
\end{equation}

This contribution represents the current redistribution induced by the magnetic coupling associated with \(i_0\). The operator
\[
R^{-1}L
\]
therefore embodies the characteristic RL time constants governing the transient evolution of the system.

Iterating the same mechanism generates the sequence
\begin{equation}
i_0,
\qquad
R^{-1}L i_0,
\qquad
(R^{-1}L)^2 i_0,
\qquad
\dots
\end{equation}

Consequently, the Krylov subspace generated by \(R^{-1}L\) defines the manifold of current distributions dynamically reachable through the coupled inductive and resistive interactions of the system, spanning the dominant space explored by the transient electromagnetic dynamics.

\section{Full solution of the null-field formulation and its MOR implementation}
\label{app:nullfield}

Starting from the coupled null-field formulation introduced in~\eqref{eq:nullfield_coupled}, the eddy-current dynamics are governed by the first equation of the coupled system, where \(\gamma\) denotes the eddy-current state, \(Y\) collects the prescribed coil currents, and \(X\) contains the control currents used to enforce the null-field condition. The objective of the problem is therefore to determine the time evolution of the control-current components
\(
X_i(t),
\ i=1,\dots,4,
\)
required to minimize the poloidal magnetic field inside the prescribed null-field region.

The transient evolution is evaluated over a temporal interval discretized with a constant time step \(\Delta t\). The corresponding full-order implicit Euler discretization is
\begin{equation}
A_K \Delta \gamma_n
=
-\Delta t\, R_K \gamma_n
+
M_Y \Delta Y_n
+
M_X \Delta X_n,
\qquad
A_K = L_K+\Delta t\,R_K ,
\label{eq:nullfield_incremental_full}
\end{equation}
with
\[
\Delta \gamma_n = \gamma_{n+1}-\gamma_n,
\qquad
\Delta Y_n = Y_{n+1}-Y_n,
\qquad
\Delta X_n = X_{n+1}-X_n.
\]

In the proposed framework, the reduced-order approximation is introduced through the MOR basis \(V_{\mathrm{r}}\):
\begin{equation}
\gamma_n \approx V_{\mathrm{r}} q_n,
\label{eq:nullfield_mor_projection}
\end{equation}
where \(q_n\) denotes the reduced-order state vector.

Substituting~\eqref{eq:nullfield_mor_projection} into~\eqref{eq:nullfield_incremental_full} and applying a Galerkin projection onto the reduced space yields
\begin{equation}
A_r \Delta q_n
=
-\Delta t\, R_r q_n
+
M_{Y,r}\Delta Y_n
+
M_{X,r}\Delta X_n,
\label{eq:nullfield_incremental_rom}
\end{equation}
with
\begin{equation}
A_r = V_{\mathrm{r}}^T A_K V_{\mathrm{r}},
\qquad
R_r = V_{\mathrm{r}}^T R_K V_{\mathrm{r}},
\label{eq:nullfield_projected_operators}
\end{equation}
and
\begin{equation}
M_{Y,r}=V_{\mathrm{r}}^T M_Y,
\qquad
M_{X,r}=V_{\mathrm{r}}^T M_X .
\label{eq:nullfield_projected_sources}
\end{equation}

The reduced update equation can therefore be written as
\begin{equation}
q_{n+1}
=
A_q q_n
+
E_q \Delta Y_n
+
F_q \Delta X_n,
\label{eq:nullfield_q_update}
\end{equation}
where
\begin{equation}
A_q =
I_r
-
\Delta t\,A_r^{-1}R_r,
\qquad
E_q = A_r^{-1}M_{Y,r},
\qquad
F_q = A_r^{-1}M_{X,r}.
\label{eq:nullfield_q_matrices}
\end{equation}

Once the MOR projection has been introduced, a further modal diagonalization can be applied at reduced-order level in order to accelerate the transient integration. Let
\begin{equation}
L_r S_r = R_r S_r \Lambda_r,
\qquad
S_r^T R_r S_r = I_r,
\qquad
S_r^T L_r S_r = \Lambda_r ,
\label{eq:nullfield_modal_eig}
\end{equation}
where
\(
\Lambda_r \in \mathbb{R}^{N_{\mathrm{mor}}\times N_{\mathrm{mor}}}
\)
is the diagonal matrix containing the reduced eigenvalues, and
\(
I_r \in \mathbb{R}^{N_{\mathrm{mor}}\times N_{\mathrm{mor}}}
\)
denotes the identity matrix. The modal reduced coordinate is defined as
\[
q_n = S_r \widetilde{q}_n .
\]

Projecting the reduced dynamics onto the modal basis yields
\begin{equation}
\widetilde{q}_{n+1}
=
(\Delta t\,I_r+\Lambda_r)^{-1}
\Lambda_r \widetilde{q}_n
+
(\Delta t\,I_r+\Lambda_r)^{-1}
\left(
\widetilde{M}_{Y,r} \Delta Y_n
+
\widetilde{M}_{X,r} \Delta X_n
\right),
\label{eq:nullfield_modal_update}
\end{equation}
where
\begin{equation}
\widetilde{M}_{Y,r} = S_r^T M_{Y,r},
\qquad
\widetilde{M}_{X,r} = S_r^T M_{X,r}.
\label{eq:nullfield_modal_sources}
\end{equation}

The modal formulation provides an equivalent accelerated implementation of the reduced update. In the following derivation, the same control equation is written in the reduced coordinates \(q_n\) for notational simplicity.

The orthogonal magnetic-field component inside the null-field region is reconstructed as
\begin{equation}
B_{\perp,n+1}
=
Q_q q_{n+1}
+
Q_Y Y_{n+1}
+
Q_X X_{n+1},
\label{eq:nullfield_B}
\end{equation}
where
\[
Q_q = Q_\gamma V_{\mathrm{r}} .
\]

Substituting the reduced update equation~\eqref{eq:nullfield_q_update} into~\eqref{eq:nullfield_B} gives
\begin{align}
B_{\perp,n+1}
&=
Q_q A_q q_n
+
Q_q E_q \Delta Y_n
+
Q_q F_q \Delta X_n
+
Q_Y Y_{n+1}
+
Q_X X_{n+1}.
\label{eq:nullfield_B_substitution}
\end{align}

Since
\[
\Delta X_n = X_{n+1}-X_n,
\]
the previous expression becomes
\begin{align}
B_{\perp,n+1}
&=
Q_q A_q q_n
+
Q_q E_q \Delta Y_n
+
Q_Y Y_{n+1}
-
Q_q F_q X_n
+
\left(
Q_X + Q_q F_q
\right)X_{n+1}.
\label{eq:nullfield_B_effective_split}
\end{align}

Defining
\begin{equation}
B^{\mathrm{eff}}_{\perp,\mathrm{known},n+1}
=
Q_q A_q q_n
+
Q_q E_q \Delta Y_n
+
Q_Y Y_{n+1}
-
Q_q F_q X_n,
\label{eq:nullfield_Bknown_eff}
\end{equation}
and
\begin{equation}
D_{\mathrm{eff}}
=
Q_X + Q_q F_q,
\label{eq:nullfield_Deff}
\end{equation}
the null-field constraint becomes
\begin{equation}
D_{\mathrm{eff}} X_{n+1}
=
-
B^{\mathrm{eff}}_{\perp,\mathrm{known},n+1}.
\label{eq:nullfield_effective_system}
\end{equation}

Finally, the control currents over the entire transient evolution are obtained through the Moore--Penrose pseudo-inverse:
\begin{equation}
X_{n+1}
=
-
D_{\mathrm{eff}}^\dagger
B^{\mathrm{eff}}_{\perp,\mathrm{known},n+1}.
\label{eq:nullfield_control_solution}
\end{equation}

Since all reduced and modal operators are precomputed offline, the resulting online stage only involves operations scaling with the reduced dimension, leading to a computational complexity proportional to \(\mathcal{O}(N_{\mathrm{mor}})\).

\section{POD Dataset Construction and Neural Reconstruction}
\label{app:POD_NN}

As anticipated, the excitation currents driving the system are described by a low-dimensional parametrization. Consequently, the transient responses also evolve within a low-dimensional subspace, allowing them to be efficiently compressed and approximated by a limited number of Proper Orthogonal Decomposition (POD) modes. The overall pipeline is split into a computationally intensive offline stage (data generation and subspace projection) and a fast online stage (neural network training and real-time inference).

The Central Solenoid (CS) and Poloidal Field (PF) excitation currents are parametrized according to a standard ramp--plateau profile:
\begin{equation}
I(t)=
\begin{cases}
A\,\dfrac{t}{t_r},
& 0 \le t < t_r,
\\
A,
& t_r \le t \le t_{\mathrm{fin}},
\end{cases}
\label{eq:ramp_plateau}
\end{equation}
where $A$ denotes the excitation amplitude, $t_r$ the rise time, and $t_{\mathrm{fin}}$ the final simulation time. This parametrization reduces the input space dimensionality to a vector $\boldsymbol{\mu}$ containing 5 geometric parameters:
\begin{equation}
\boldsymbol{\mu} = [A_{cs},\, t_{r,cs},\, A_{pf},\, t_{r,pf},\, t_{\mathrm{fin}}]^T \in \mathbb{R}^5.
\end{equation}

\begin{figure}[H]
    \centering
    \includegraphics[width=0.62\textwidth]{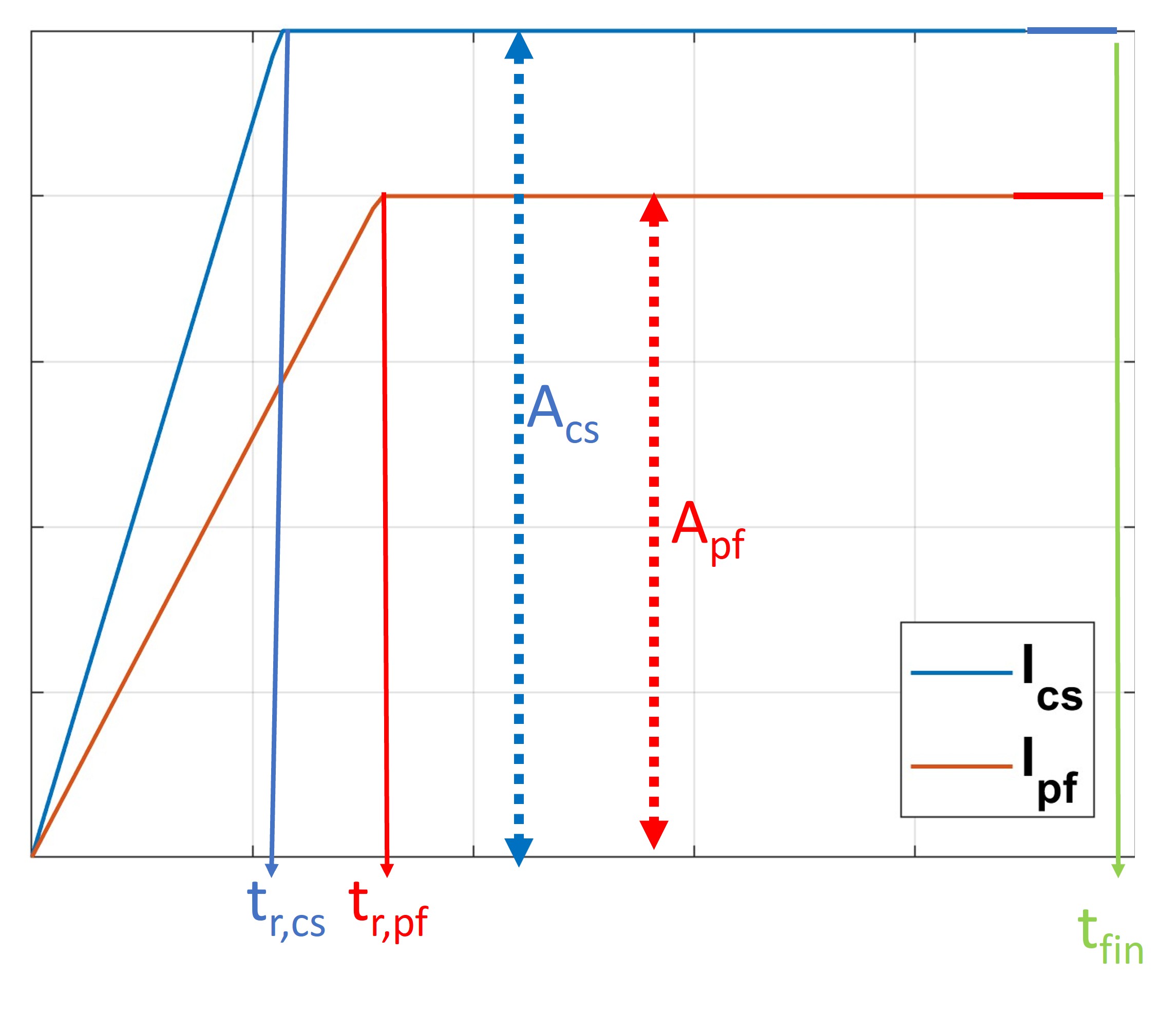}
    \caption{Ramp--plateau parametrization of the excitation currents and corresponding definition of the parameter vector \(\mu\), used as input of the NN.}
    \label{fig:ramp-plateau}
\end{figure}

The five parameters defining the NN input vector\(\boldsymbol{\mu}\) and the corresponding ramp--plateau excitation profiles are illustrated in Figure~\ref{fig:ramp-plateau}.

For a given parameter configuration $\boldsymbol{\mu}$, the continuous transient electromagnetic response $\mathbf{X}(t;\boldsymbol{\mu})$ encompassing the 4 control currents can be approximated via POD expansion as:
\begin{equation}
\mathbf{X}(t;\boldsymbol{\mu}) \approx \mathbf{X}_{mean}(t) + \sum_{j=1}^{N_{\mathrm{POD}}} a_i(\boldsymbol{\mu})\,\boldsymbol{\psi}_i(t),
\end{equation}
where $\mathbf{X}_{mean}(t)$ is the mean transient trajectory. The continuous modal functions $\boldsymbol{\psi}_i(t)$ form an orthonormal basis in the time domain, satisfying the condition:
\begin{equation}
\label{eq:pod_orthogonality}
\int_{0}^{t_{\mathrm{fin}}} \boldsymbol{\psi}_i(t) \cdot \boldsymbol{\psi}_j(t) \, \mathrm{d}t = \delta_{ij}
\end{equation}
where $\delta_{ij}$ is the Kronecker delta. The column vector $\mathbf{a}(\boldsymbol{\mu}) = [a_1(\boldsymbol{\mu}), \dots, a_r(\boldsymbol{\mu})]^T \in \mathbb{R}^r$ represents the reduced coordinates in the POD subspace of dimension \(N_{\mathrm{POD}}\).

To handle the data numerically, a uniform temporal grid with $n_t$ time steps is introduced:
\begin{equation}
t_i = \Delta t \, (i-1), \qquad i=1,\dots,n_t.
\end{equation}
For each realization $\boldsymbol{\mu}^{(k)}$, the 4 transient responses sampled over the grid are concatenated into a single snapshot vector:
\begin{equation}
\mathbf{x}_c(\boldsymbol{\mu}^{(k)}) = [X_1(t_1), \dots, X_1(t_{n_t}), X_2(t_1), \dots, X_4(t_{n_t})]^T \in \mathbb{R}^{4n_t}.
\end{equation}

The discrete POD reconstruction of a snapshot is written as
\begin{equation}
\mathbf{x}_c
=
\mathbf{\Psi}_{\mathrm{mat}}\,\mathbf{a},
\qquad
\mathbf{\Psi}_{\mathrm{mat}}
=
[\boldsymbol{\psi}_1,\dots,\boldsymbol{\psi}_{N_{\mathrm{POD}}}]
\in
\mathbb{R}^{4n_t\times N_{\mathrm{POD}}},
\end{equation}
where, from the orthogonality condition~\eqref{eq:pod_orthogonality},
\begin{equation}
\mathbf{\Psi}_{\mathrm{mat}}^T
\mathbf{\Psi}_{\mathrm{mat}}
=
\mathbf{I}_{N_{\mathrm{POD}}},
\end{equation}
with
\(
\mathbf{I}_{N_{\mathrm{POD}}}
\in
\mathbb{R}^{N_{\mathrm{POD}}\times N_{\mathrm{POD}}}
\)
denoting the identity matrix.
By collecting $N_s$ transient realizations from the MOR solver, a global snapshot matrix is built:
\begin{equation}
\mathbf{X}_{mat} = [\mathbf{x}_c(\boldsymbol{\mu}^{(1)}), \, \dots, \, \mathbf{x}_c(\boldsymbol{\mu}^{(N_s)})] \in \mathbb{R}^{4n_t \times N_s}.
\end{equation}
After subtracting the ensemble mean matrix to center the dataset ($\widetilde{\mathbf{X}} = \mathbf{X}_{mat} - \mathbf{x}_{mean}\mathbf{1}_{N_s}^T$), the dominant modes and coefficients are extracted. 

This compression can be mathematically formalized via Singular Value Decomposition (SVD),
\[
\widetilde{\mathbf{X}}
=
\mathbf{U}\mathbf{\Sigma}\mathbf{V}^T,
\]
or equivalently via a truncated Column-Pivoted QR decomposition,
\[
\widetilde{\mathbf{X}}
=
\mathbf{Q}\mathbf{R}.
\]
The POD basis matrix is then identified as
\[
\mathbf{\Psi}_{\mathrm{mat}}
\equiv
\mathbf{U}_{N_{\mathrm{POD}}}
\equiv
\mathbf{Q}_{N_{\mathrm{POD}}},
\]
where
\(
\mathbf{U}_{N_{\mathrm{POD}}}
\)
and
\(
\mathbf{Q}_{N_{\mathrm{POD}}}
\)
contain the first \(N_{\mathrm{POD}}\) dominant modes. The reduced training coefficients matrix \(\mathbf{A}_{\mathrm{mat}}\) is extracted by projecting the centered snapshots onto the orthonormal POD basis:
\begin{equation}
\mathbf{A}_{\mathrm{mat}}
=
\mathbf{\Psi}_{\mathrm{mat}}^T
\widetilde{\mathbf{X}}
=
[
\mathbf{a}(\boldsymbol{\mu}^{(1)}),
\dots,
\mathbf{a}(\boldsymbol{\mu}^{(N_s)})
]
\in
\mathbb{R}^{N_{\mathrm{POD}}\times N_s}.
\end{equation}
A supervised Artificial Neural Network (NN) is built to learn the low-dimensional mapping $\boldsymbol{\mu} \mapsto \mathbf{a}(\boldsymbol{\mu})$. Before training, both inputs and targets are scaled using standard normalization:
\begin{equation}
\boldsymbol{\mu}_n = \frac{\boldsymbol{\mu} - \mathbf{m}_{\mu}}{\boldsymbol{\sigma}_{\mu}}, \qquad \mathbf{a}_n = \frac{\mathbf{a} - \mathbf{m}_{a}}{\boldsymbol{\sigma}_{a}}.
\end{equation}
The optimal NN parameters are obtained by minimizing the Mean Squared Error (MSE) loss function over the training dataset $\mathbf{A}_{mat}$.

Training the neural network to predict the reduced vector \(\mathbf{a} \in \mathbb{R}^{N_{\mathrm{POD}}}\),rather than the full transient vector $\mathbf{x}_c \in \mathbb{R}^{4n_t}$ yields major algorithmic advantages. Since the POD dimension \(N_{\mathrm{POD}}\) is very small (typically on the order of ten) compared to the physical discretization size ($4n_t \sim 10^3-10^4$), the output layer dimension is significantly reduced. This low-dimensional target space prevents overfitting, ensures smooth loss landscapes, and drives the neural network to full convergence with minimal training epochs.

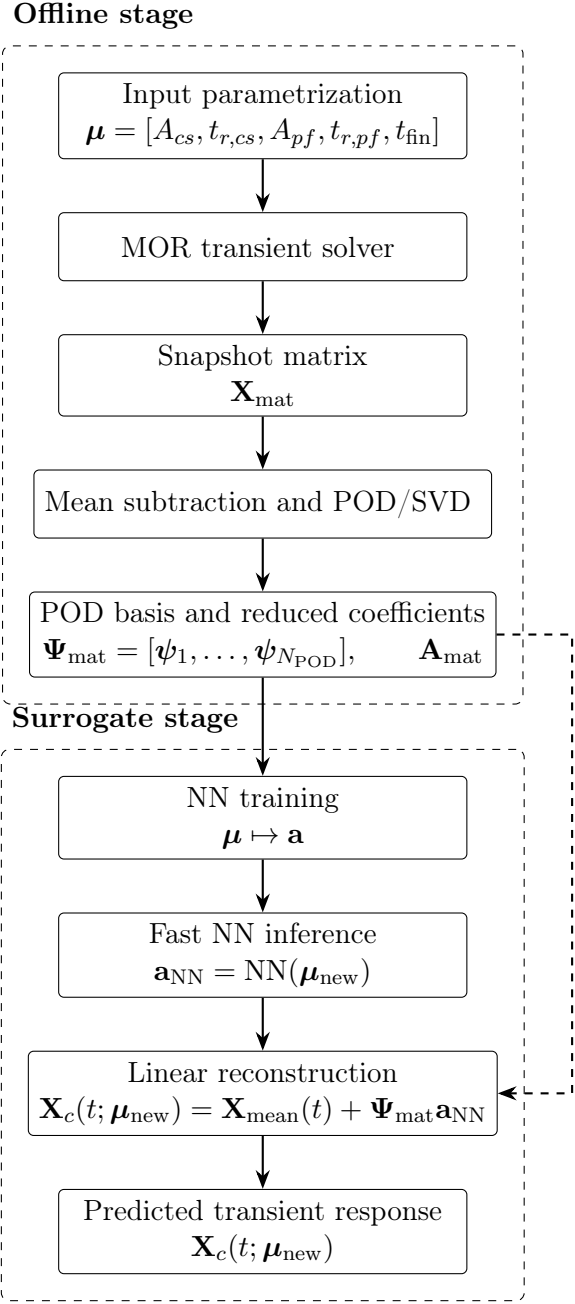
\begin{figure}[H]
\centering
\begin{tikzpicture}[
    node distance=0.7cm,
    block/.style={
        draw,
        rectangle,
        rounded corners=2pt,
        minimum width=5.4cm,
        minimum height=0.9cm,
        align=center,
        font=\small
    },
    arrow/.style={-{Stealth[length=2.5mm]}, thick},
    dashedarrow/.style={-{Stealth[length=2.5mm]}, thick, dashed},
    group/.style={
        draw,
        dashed,
        rounded corners=3pt,
        inner sep=0.35cm
    }
]

\node[block] (mu) {
Input parametrization
\\
\(
\boldsymbol{\mu}
=
[A_{cs},t_{r,cs},A_{pf},t_{r,pf},t_{\mathrm{fin}}]
\)
};

\node[block, below=of mu] (mor) {
MOR transient solver
};

\node[block, below=of mor] (snap) {
Snapshot matrix
\\
\(
\mathbf{X}_{\mathrm{mat}}
\)
};

\node[block, below=of snap] (pod) {
Mean subtraction and POD/SVD
};

\node[block, below=of pod] (basis) {
POD basis and reduced coefficients
\\
\(
\mathbf{\Psi}_{\mathrm{mat}}
=
[\boldsymbol{\psi}_1,\ldots,\boldsymbol{\psi}_{N_{\mathrm{POD}}}],
\qquad
\mathbf{A}_{\mathrm{mat}}
\)
};

\draw[arrow] (mu) -- (mor);
\draw[arrow] (mor) -- (snap);
\draw[arrow] (snap) -- (pod);
\draw[arrow] (pod) -- (basis);

\node[
    group,
    fit=(mu) (mor) (snap) (pod) (basis)
] (offlinebox) {};

\node[
    anchor=south west,
    font=\bfseries\small,
    yshift=2pt
]
at (offlinebox.north west)
{Offline stage};

\node[block, below=1.3cm of basis] (train) {
NN training
\\
\(
\boldsymbol{\mu}
\mapsto
\mathbf{a}
\)
};

\node[block, below=of train] (infer) {
Fast NN inference
\\
\(
\mathbf{a}_{\mathrm{NN}}
=
\mathrm{NN}(\boldsymbol{\mu}_{\mathrm{new}})
\)
};

\node[block, below=of infer] (rec) {
Linear reconstruction
\\
\(
\mathbf{X}_c(t;\boldsymbol{\mu}_{\mathrm{new}})
=
\mathbf{X}_{\mathrm{mean}}(t)
+
\mathbf{\Psi}_{\mathrm{mat}}
\mathbf{a}_{\mathrm{NN}}
\)
};

\node[block, below=of rec] (output) {
Predicted transient response
\\
\(
\mathbf{X}_c(t;\boldsymbol{\mu}_{\mathrm{new}})
\)
};

\draw[arrow] (basis) -- (train);
\draw[arrow] (train) -- (infer);
\draw[arrow] (infer) -- (rec);
\draw[arrow] (rec) -- (output);

\draw[dashedarrow]
(basis.east) -- ++(1.0cm,0)
|- (rec.east);

\node[
    group,
    fit=(train) (infer) (rec) (output)
] (onlinebox) {};

\node[
    anchor=south west,
    font=\bfseries\small,
    yshift=2pt
]
at (onlinebox.north west)
{Surrogate stage};

\end{tikzpicture}

\caption{
Flowchart of the proposed POD--NN framework.
The offline stage generates the snapshot database and extracts the POD basis and reduced coefficients through SVD compression.
The surrogate stage learns the parametric mapping between the excitation parameters and the reduced coordinates, enabling fast transient reconstruction.
}
\label{fig:app_flowchart_pod_nn}
\end{figure}

During the online operational stage, the evaluation of a new parameter configuration $\boldsymbol{\mu}_{new}$ bypasses the expensive full-order or reduced-order numerical solvers. The process is instantaneous and structured as follows: the network computes the single reduced column vector $\mathbf{a}_{nn} = \text{NN}(\boldsymbol{\mu}_{new})$, and the high-dimensional physical transient response is reconstructed via a single matrix-vector multiplication and mean re-addition:
\begin{equation}
\mathbf{x}_c^{pred} = \mathbf{x}_{mean} + \mathbf{\Psi}_{mat} \mathbf{a}_{nn}.
\end{equation}
The entire POD--NN workflow, including both the offline dataset construction and the online surrogate inference stages, is schematically summarized in the flowchart reported in Figure~\ref{fig:app_flowchart_pod_nn}.

\section*{Acknowledgements}
The author gratefully acknowledges Pierre Bauer for providing the numerical data and technical information concerning the ITER feeder busbar magnetic-field studies.

The author also dedicates this work to his sons: Ada, Stefano, and Francesco.



\end{document}